%% file: param-abstract.tex
\documentclass[runningheads,a4paper,envcountsame,envcountsect]{llncs}

\usepackage{amssymb}
\usepackage{graphicx}

\usepackage{url}

\let\saveendproof=\endproof
\def\endproof{\qed\saveendproof}

\newtheorem{obs}[theorem]{Observation}

\usepackage{booktabs}
\usepackage{amssymb}
\usepackage{graphicx}
\usepackage{epstopdf}
\usepackage{url}
\usepackage{color}
\usepackage[usenames,dvipsnames,svgnames,table]{xcolor}

\usepackage{tikz}

\usepackage{pgfplots} 

\usetikzlibrary{positioning}

\def\visible<#1>{}  

\usepackage[utf8]{inputenc}

\usepackage[english]{babel}
\usepackage{amsfonts}
\usepackage{amsmath}
\usepackage{latexsym}
\usepackage{color}
\usepackage{subfigure}
\usepackage{enumerate}

\usepackage{hyperref}  

\usepackage{ifpdf}

\graphicspath{{figures/}}

\newcommand{\old}[1]{{}}

\newcommand{\bb}{\mathbb}

\newcommand{\R}{\bb R}
\newcommand{\Q}{\bb Q}
\newcommand{\Z}{\bb Z}

\def\ve#1{\mathchoice{\mbox{\boldmath$\displaystyle\bf#1$}}
{\mbox{\boldmath$\textstyle\bf#1$}}
{\mbox{\boldmath$\scriptstyle\bf#1$}}
{\mbox{\boldmath$\scriptscriptstyle\bf#1$}}}

\newenvironment{psmallmatrixbig}{\bigl(\smallmatrix}{\endsmallmatrix\bigr)}
\newcommand\InlineFrac[2]{#1/#2}  
\newcommand\ColVec[3][\relax]
{
  \ifx#1\relax
  \bgroup\let\frac=\InlineFrac\begin{psmallmatrixbig}#2\vphantom{/}\\#3\vphantom{/}\end{psmallmatrixbig}\egroup
  \else
  \bgroup\let\frac=\InlineFrac\begin{psmallmatrixbig}\ifx#200\else#2/#1\fi\\\ifx#300\else#3/#1\fi\end{psmallmatrixbig}\egroup
  \fi
}

\makeatletter
\renewcommand{\pod}[1]
{\allowbreak\mathchoice{\mkern18mu}{\mkern8mu}{\mkern8mu}{\mkern8mu}(#1)}
\makeatother

\chardef\Myunderscore=`\_
\newcommand\underscore{\Myunderscore\allowbreak}
\let\_=\underscore

\newcommand\githubsearchurl{https://github.com/mkoeppe/infinite-group-relaxation-code/search}
\input{sage-commands}
\DeclareRobustCommand\sage[1]{\texttt{#1}}
\DeclareRobustCommand\sagefunc[1]{\pgfkeys{/sagefunc/#1}}

\usepackage[normalem]{ulem}

\usepackage[]{algorithm2e}
\usepackage{booktabs,array,ragged2e}
\usepackage{pdflscape}
\usepackage{blkarray}
\definecolor{mediumspringgreen}{rgb}{0.0, 0.98039215, 0.60392156}

\def\groupcodedir{.}

\newcommand\SymbolicElement[2]{{#1}\,|_{={#2}}}

\title{Toward computer-assisted discovery and automated proofs of cutting plane theorems}
\titlerunning{Discovery and proof of cutting plane theorems}

\author{Matthias K\"oppe \and
  Yuan Zhou}

\institute{Dept.\ of Mathematics, University of California, Davis\\
  \texttt{mkoeppe@math.ucdavis.edu}, \texttt{yzh@math.ucdavis.edu}}

\date{\today}

\usepackage[normalem]{ulem}

\begin{document}

\maketitle
\begin{abstract}
  Using a metaprogramming technique and semialgebraic computations, we provide
  computer-based proofs for old and new cutting-plane theorems in 
  Gomory--Johnson's model of cut generating functions.
\end{abstract}

\section{Introduction}

Inspired by the spectacular breakthroughs of the polyhedral method for
combinatorial optimization in the 1980s, generations of researchers have
studied the facet structure of convex hulls to develop strong cutting planes.
It is a showcase of the power of experimental mathematics: Small examples are
generated, their convex hulls are computed (for example, using the popular
tool PORTA \cite{porta}), conjectures are formed, theorems are proved.  Some proofs feature
brilliant new ideas; other proofs are routine.  Once the theorems have been
found and proved, separation algorithms for the cutting planes are
implemented.  Numerical tests are run, the strength-versus-speed trade-off is
investigated, parameters are tuned, papers are written.

In this paper, we ask how much of this process can be automated:
In particular, \emph{can we use algorithms to discover and prove theorems about
cutting planes?}  
This paper is part of a larger project in which we aim to automate more stages
of this pipeline.  
We focus on general integer and mixed integer programming, rather than
combinatorial optimization, and use the framework of cut-generating functions
\cite{conforti2013cut}, specifically those of the classic single-row Gomory--Johnson
model \cite{infinite,infinite2}.  Cut-generating functions are an attractive
framework for our study for several reasons.  First, it is essentially
dimensionless: Cuts obtained from cut-generating functions can be applied to
problems of arbitrary dimension.  Second, it may be a way to the mythical
multi-row cuts, sometimes dubbed the ``holy grail of integer programming,''
though the computational approaches so far have disappointed.  Third, work on
new cuts in the single-row Gomory--Johnson model has, with few exceptions
, become a
routine, but error-prone task that leads to proofs of enormous complexity; see
for example \cite{Miller-Li-Richard2008,Richard-Li-Miller-2009:Approximate-Liftings}.  Fourth, finding new
cuts in the multi-row Gomory--Johnson model has a daunting complexity, and few
attempts at a systematic study have been made.  Fifth, working on the
Gomory--Johnson model is timely because only recently, after decades of
theoretical investigations, the first
computational tools for cut-generating functions in this model became
available in \cite{basu-hildebrand-koeppe:equivariant} and
the software implementation \cite{infinite-group-relaxation-code}.

Of course, automated theorem proving is not a new proposition.  Probably the
best known examples in the optimization community are the proof of the Four
Color Theorem, by Appel--Haken \cite{appel1977}
, and more
recently and most spectacularly the proof of the Kepler Conjecture by Hales
\cite{hales2005kepler} and again within Hales' Flyspeck project in
\cite{hales2015formal}. 
In the domains of combinatorics, number theory, and plane geometry, Zeilberger
with long-term collaborator Shalosh B.~Ekhad have pioneered automated
discovery and proof of theorems; see, for example \cite{zeilberger-gt}.  Many
sophisticated automated theorem provers, by names 
such are HOL light, Coq, Isabelle, Mizar, etc. are available nowadays; see
\cite{freek-comparison} and the references within for an interesting overview.

Our approach is pragmatic.  Our theorems and proofs come from a
metaprogramming trick, applied to the practical software implementation
\cite{infinite-group-relaxation-code} of computations with the Gomory--Johnson
model; followed by computations with semialgebraic cell complexes.  As such,
all of our techniques are reasonably close to mathematical programming
practice.  The correctness of all of our proofs depends on the correctness of
the underlying implementation.  We make no claims that our proofs can be
formalized in the sense of the above mentioned formal proof systems that break
every theorem down to the axioms of mathematics; in fact, we make no attempt
to even use an automated theorem proving system.

Our software is in an early, proof-of-concept stage of development.  In this
largely computational and experimental paper we report on the early successes
of the software.  We computationally verify the results on the
\sagefunc{gj_forward_3_slope}\footnote{A function name shown
  in typewriter font is the name of the constructor of this function in the
  Electronic Compendium, part of the SageMath
  program~\cite{infinite-group-relaxation-code}.
  In an online copy of this paper, there are hyperlinks that lead to a search
  for this function in the GitHub repository.}
and \sagefunc{drlm_backward_3_slope}
functions.
We find a correction to a theorem by Chen \cite{chen} regarding the
extremality conditions for his \sagefunc{chen_4_slope} family.\footnote{This
  is a new result, which should not be confused with our previous result in
  \cite{zhou:extreme-notes} regarding Chen's family of 3-slope functions
  (\sagefunc{chen_3_slope_not_extreme}).}  
We find a correction to a result by Miller, Li and Richard \cite{Miller-Li-Richard2008} on
the so-called $\mathrm{CPL}_3^=$ functions (\sage{mlr\_cpl3\_}\dots). 
We discover several new parametric
families, \sagefunc{kzh_3_slope_param_extreme_1} and
\sagefunc{kzh_3_slope_param_extreme_2}, of extreme functions and
corresponding theorems regarding their extremality, with automatic proofs.

These new theorems are entirely unremarkable; the plan for the future is to
make up for it by sheer quantity.\footnote{With this anticlimactic sentence,
  the introduction ends.  We leave it to the reader to speculate about the
  possible computational implications of having a large library of families of
  cut-generating functions available.  For example, could the diversity of the
  cuts offered by such a library, and their rich parametrization, become
  crucial to the success of a branch-and-cut algorithm auto-tuned by machine
  learning?
}

\clearpage
\section{The Gomory--Johnson model}

We restrict ourselves to the single-row (or, ``one-dimensional'') infinite
group problem, which has attracted most of the attention in the past and for
which the software~\cite{infinite-group-relaxation-code} is available.  It can
be written as
\begin{equation}
  \label{GP} 
  \begin{aligned}
    &\sum_{r \in \R} r\, y(r) \equiv f \pmod{1}, \\
    &y\colon \R\to\Z_+ \text{ is a function of finite support}, 
  \end{aligned}
\end{equation}
where $f$ is a given element of $\R\setminus \Z$. 
We study the convex hull $R_{f}(\R,\Z)$ of
the set of all functions $y\colon \R \to \Z_+$ satisfying the constraints
in~\eqref{GP}.  The elements of the convex hull are understood as functions
$y\colon \R \to \R_+$. 

After a normalization, valid inequalities for the convex set $R_{f}(\R,\Z)$
can be described using so-called \emph{valid functions} $\pi\colon \R\to\R$
via $\langle \pi, y \rangle := \sum_{r \in \R} \pi(r)y(r) \geq 1$.  
Valid functions~$\pi$ are cut-generating functions for pure integer programs.
Take a row of the optimal simplex tableau of an integer program, corresponding
to a basic variable~$x_i$ that currently takes a fractional value: 
$$ x_i = -f_i + \sum_{j\in N} r_j x_j, \quad x_i \in \Z_+,\ \ve x_N \in
\Z_+^N.  $$
Then a valid function $\pi$ for  $R_{f_i}(\R,\Z)$ gives a valid inequality 
$ \sum_{j\in N} \pi(r_j) x_j \geq 1 $
for the integer program. 
(By a theorem of Johnson \cite{johnson}, this extends easily to the mixed integer case:
A function $\psi$ can be associated to $\pi$, so that they together form a
\emph{cut-generating function pair} $(\psi, \pi)$, which gives the
coefficients of the continuous and of the integer variables.)

In the finite-dimensional case, instead of merely valid inequalities, one is
interested in stronger inequalities such as tight valid inequalities and
facet-defining inequalities.  These r\^oles are taken in our
infinite-dimensional setting by \emph{minimal functions} and \emph{extreme
  functions}.  Minimal functions are those valid functions that are pointwise
minimal; extreme functions are those that are not a proper convex combination
of other valid functions.

By a theorem of Gomory and
Johnson~\cite{infinite}, minimal functions for $R_f(\R,\Z)$ are classified:  
They are the subadditive functions $\pi\colon \R\to\R_+$ that are
periodic modulo~$1$ and satisfy the \emph{symmetry condition} $\pi(x) +
\pi(f - x) = 1$ for all $x\in\R$.  

Obtaining a full classification of the \emph{extreme} functions has proved to
be elusive, however various authors have defined parametric families of
extreme functions and provided extremality proofs for these families.  These
parametric families of extreme functions from the literature, as well as
``sporadic'' extreme functions, have been collected in an electronic
compendium as a part of the software~\cite{infinite-group-relaxation-code}; see \cite{zhou:extreme-notes}.

We refer the interested reader to the recent surveys
\cite{corner_survey,igp_survey,igp_survey_part_2} for a more
detailed exposition. 

\section{Examples of cutting-plane theorems in the Gomory--Johnson model}

To illustrate what cutting-plane theorems in the Gomory--Johnson model look
like, we give three examples of such theorems, paraphrased for precision from
the literature where they were stated.  As it turns out, the last theorem is incorrect.

\setlength{\belowcaptionskip}{-10pt}

\begin{figure}[h]
\centering
\includegraphics[width=\linewidth]{\groupcodedir/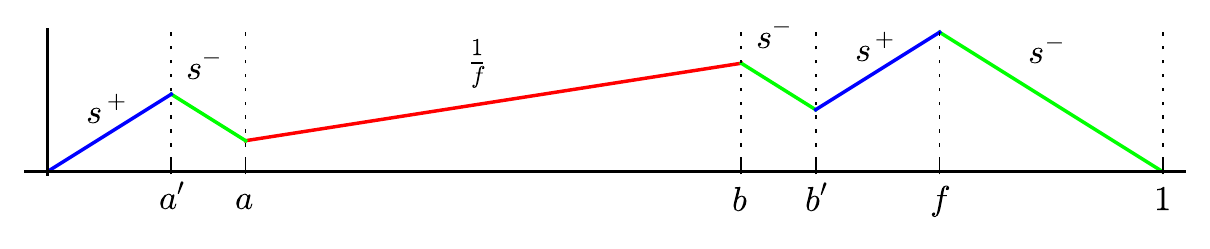}\par\vspace{-3ex}
\caption{\sagefunc{gj_forward_3_slope}}
\end{figure}
\begin{theorem}[reworded from Gomory--Johnson {\cite[Theorem 8]{tspace}}]
  \label{th:gj_forward_3_slope}
Let $f \in (0,1)$ and $\lambda_1, \lambda_2 \in \R$. Define the periodic,
piecewise linear \sagefunc{gj_forward_3_slope} function $\pi\colon \R/\Z
\to \R$ as follows.   
The function $\pi$ satisfies $\pi(0) = \pi(1) =0$; it has $6$ pieces between $0$ and $1$ with breakpoints at $0, a', a, b, b', f$ and $1$, where 
$a = \frac{\lambda_1 f}{2},  \; a' = a + \frac{\lambda_2 (f-1)}{2}, \; b = f - a$ and $ b' = f - a'$.
The slope values of $\pi$ on these pieces are $s^+, s^-, \frac{1}{f}, s^-, s^+$ and $s^-$, repectively, where $s^+= \frac{\lambda_1+\lambda_2}{\lambda_1 f + \lambda_2 (f-1)}$ and $s^- = \frac{1}{f-1}$.
If the parameters $\lambda_1$ and $\lambda_2$ satisfy that \emph{(i)} $0 \leq \lambda_1
\leq \frac{1}{2}$, \emph{(ii)} $0 \leq \lambda_2 \leq 1$
and 
\emph{(iii)} 
$0 < \lambda_1 f + \lambda_2 (f-1) $
, then the function $\pi$ is an
extreme function for $R_f(\R/\Z)$.
\end{theorem}


\begin{figure}[h]
\centering
\includegraphics[width=\linewidth]{\groupcodedir/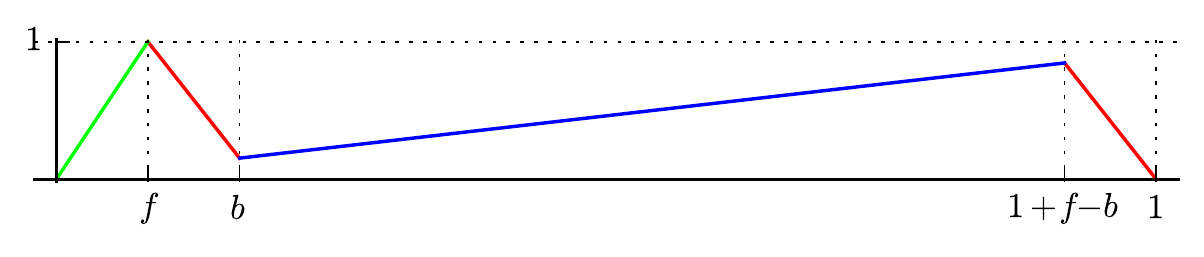}\par\vspace{-3ex}
\caption{\sagefunc{drlm_backward_3_slope}}
\end{figure}
\begin{theorem}[Dey--Richard--Li--Miller \cite{dey1}; in this form, for the
  real case, in~{\cite[Theorem 4.1]{zhou:extreme-notes}}]
  \label{th:drlm_backward_3_slope}
Let $f$ and $b$ be real numbers such that $0 < f < b \leq \frac{1+f}{4}$. The periodic,
piecewise linear \sagefunc{drlm_backward_3_slope} function $\pi\colon \R/\Z\to \R$ defined as
follows is an extreme function for $R_f(\R/\Z)$:
\[
\pi(x) =
  \begin{cases}
   \frac{x}{f} & \text{if } 0 \leq x \leq f \\
   1 + \frac{(1+f-b)(x-f)}{(1+f)(f-b)} & \text{if } f \leq x \leq b \\
   \frac{x}{1+f} & \text{if } b \leq x \leq 1+f-b \\
   \frac{(1+f-b)(x-1)}{(1+f)(f-b)}  & \text{if } 1+f-b \leq x \leq 1
  \end{cases}
\]
\end{theorem}

\begin{figure}[h]
\centering
\includegraphics[width=\linewidth]{\groupcodedir/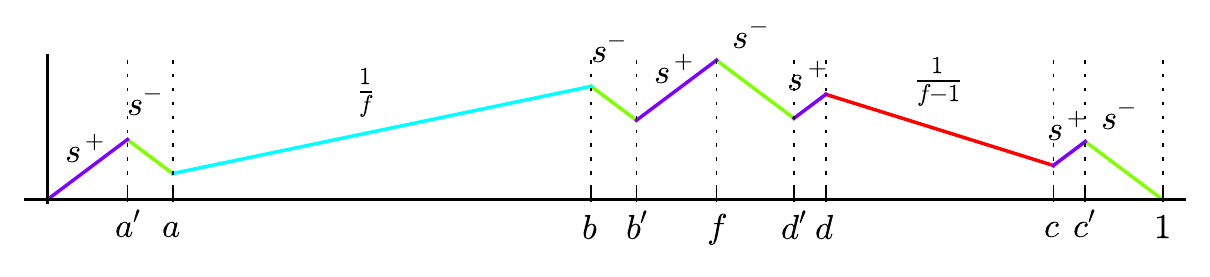}\par\vspace{-3ex}
\caption{\sagefunc{chen_4_slope}}
\end{figure}
\begin{theorem}[reworded from Chen {\cite[Theorem 2.2.1]{chen}}]
  \label{th:chen_4_slope}
Let $f \in (0,1)$, $s^+ >0,  s^-  < 0$ and $\lambda_1, \lambda_2 \in
\R$. Define the periodic, piecewise linear \sagefunc{chen_4_slope} function $\pi\colon \R/\Z \to \R$ as follows.  
The function $\pi$ satisfies $\pi(0) = \pi(1) =0$; it has $10$ pieces between $0$ and $1$ with breakpoints at $0, a', a, b, b', f, d', d, c, c',1$, where 
\[a' = \frac{\lambda_1 (1- s^- f)}{2(s^+ -s^-)}, \quad
a = \frac{\lambda_1 f}{2}, \quad
c = 1 - \frac{\lambda_2 (1 -f)}{2}, \quad
c'  = 1 -  \frac{\lambda_2 (1 - s^+(1-f))}{2(s^+ - s^-)}\]
and  
\(b = f - a\),  \( b' = f - a'\) ,  \( d = 1 + f - c \), \( d'  = 1 + f - c'\).
The slope values of $\pi$ on these pieces are $s^+, s^-, \frac{1}{f}, s^-, s^+, s^-, s^+, \frac{1}{f-1}, s^+$ and $s^-$, repectively.
If the parameters $f, \lambda_1, \lambda_2, s^+$ and  $s^-$ satisfy that
\[ f \geq \frac{1}{2}, \quad
 s^+ \geq \frac{1}{f}, \quad
 s^-  \leq \frac{1}{f-1},\quad
 0 \leq \lambda_1 < \min \lbrace \frac{1}{2}, \frac{s^+ - s^-}{s^+(1-s^- f)}\rbrace, \text{ and}\]
\[ f - \frac{1}{s^+} < \lambda_2 < \min\Big\lbrace \frac{1}{2}, \frac{s^+ -s^-}{s^-(s^+(f-1)-1)} \Big\rbrace,\]
then $\pi$ is an extreme function for $R_f(\R,\Z)$.
\end{theorem}



\begin{obs}
  \begin{enumerate}[(i)]
  \item These theorems are about families of periodic, continuous piecewise
    linear functions $\pi\colon \R\to\R$ that depend on a finite number of
    real parameters in a way that breakpoints and slope values can be written
    as rational functions of the parameters.
  \item 
    There are natural conditions on the parameters to make the function even
    constructible; for example, in \autoref{th:drlm_backward_3_slope}, if
    $f < b$ is violated, then the function is not well-defined.  These
    conditions are inequalities of rational functions of the parameters.
    Hence the set of parameter tuples such that the construction describes a function is a \emph{semialgebraic set}.
  \item 
    There are additional conditions on the parameters that ensure that the
    function is an extreme function.  Again, all of these conditions are
    inequalities of rational functions of the parameters.
    Hence the set of parameter tuples such that the construction gives an extreme function is a \emph{semialgebraic set}.
  \end{enumerate}
\end{obs}

\begin{remark}
  Some families of extreme functions in the literature are defined in more
  general ways. Some use parameters that are integers (for example,
  \sagefunc{drlm_2_slope_limit} has integer parameters that control the number
  of pieces of the function). Others use non-algebraic operations such as 
  the floor/ceiling/fractional part operations to define the
  breakpoints and slope values of the function (for example,
  \sagefunc{dg_2_step_mir}).  Another family, \sagefunc{bhk_irrational},
  requires an arithmetic condition, the $\Q$-linear independence of certain
  parameters, for extremality.
  These families are beyond the scope of this paper.
\end{remark}

\section{Semialgebraic cell structure of extremality proofs}
\label{s:semialgebraic-cells}

The minimality of a given periodic piecewise linear function can be easily
tested algorithmically; see, for example, \cite[Theorem 3.11]{igp_survey}.  
Basu, Hildebrand, and K\"oppe \cite{basu-hildebrand-koeppe:equivariant}
gave the first algorithmic tests for extremality for a given function~$\pi$ whose
breakpoints are rational with a common denominator~$q$. 
The simplest of these tests uses their finite-oversampling theorem (see
\cite[Theorem 8.6]{igp_survey_part_2} for its strongest form). 
Extremality of the function~$\pi$ is equivalent to the extremality of its
restriction to the refined grid $\frac1{3q}\Z/\Z$ for the finite master group
problem. Thus it can be tested by finite-dimensional linear algebra.

The proof of the finite-oversampling theorem in
\cite{basu-hildebrand-koeppe:equivariant} (see also \cite[section
7.1]{igp_survey_part_2} for a more high-level exposition) provides another
algorithm, based on the computation of ``affine-imposing'' (``covered'')
intervals and the construction of ``equivariant'' perturbation functions.
This algorithm in \cite{basu-hildebrand-koeppe:equivariant} is also tied to
the use of the grid $\frac1q\Z/\Z$; but it has since been generalized in the
practical implementation \cite{infinite-group-relaxation-code} to give a
completely \emph{grid-free algorithm}, which is suitable also for rational breakpoints
with huge denominators and for irrational breakpoints.\footnote{The finiteness
  proof of the algorithm, however, does depend on the rationality of the
  data. In this paper we shall ignore the case of functions with non-covered intervals
  and irrational breakpoints, such as the \sagefunc{bhk_irrational} family,
  which necessitates testing the $\Q$-linear
  independence of certain parameters.}

\begin{obs}
  \label{ob:algebraic-algorithm}
  On close inspection of this grid-free algorithm, we see that it only uses
  algebraic operations, comparisons, and branches
  (\emph{if}-\emph{then}-\emph{else} statements and loops), and then returns
  either \emph{True} (to indicate extremality) or \emph{False}
  (non-extremality).
\end{obs}

Enter \emph{parametric analysis} of the algorithm, that is, we wish to run the
algorithm for a function from a parametric family and observe how the run of
the algorithm and its answer changes, depending on the parameters.  It is then
a simple observation that for any algorithm of the type described in
Observation~\ref{ob:algebraic-algorithm}, the set of parameters where the algorithm
returns \emph{True} must be a union of sets described by equations
and inequalities of rational functions in the parameters. If the number of operations (and thus the number of
branches) that the algorithm takes is bounded finitely, then it will be a
finite union of ``cells'', each corresponding to a particular outcome of
comparisons that led to branches, and each described by finitely many equations and
inequalities of rational functions in the parameters.  Thus it will be a semialgebraic
set.  

Within each of the cells, we get the ``same'' proof of extremality.  A
complete proof of extremality for a parametric family is merely a collection
of cells, with one proof for each of them.  This is what we compute as we
describe below.

\section{Computing one proof cell by metaprogramming}
\label{s:one-cell-metaprogramming}

Now we describe how we compute one cell of the proof.  We assume that we are
given a tuple of \emph{concrete} parameter values; we will compute a semialgebraic
description of a cell of parameter tuples for which the algorithm takes the
same branches.

It is well known that modern programming languages and systems provide
facilities known as ``operator overloading'' or ``virtual methods'' that allow
us to conveniently write ``generic'' programs that can be run on objects of
various types. For example, the program \cite{infinite-group-relaxation-code},
which is written in the SageMath system \cite{sage}, by default works with
(arbitrary-precision) rational numbers; but when parameters are irrational
algebraic numbers, it constructs a suitable number field, embedded into the
real numbers, and makes exact computations with the elements of this number
field.\footnote{The details of this are not relevant to the present paper, 
  but one reader of a version of this paper was intrigued by this sentence, so 
  we elaborate.  These number fields are algebraic field extensions (of some
  degree $d$) of the field~$\Q$ of rational numbers, in much the same way that the
  field~$\mathbb C$ of complex numbers is an algebraic field extension (of degree $d=2$) of
  the field~$\R$ of real numbers.  Elements of the field are represented as a
  coordinate vector of dimension~$d$ over the base field, and all arithmetic computations
  are done by manipulating these vectors.  The computational overhead,
  compared to arithmetic over~$\Q$, is negligible for quadratic field
  extensions thanks to a specialized implementation in SageMath, and 
  within a factor of $10$ at least for field extensions of degree $d \leq 60$, as tested
  by \sage{\%timeit \sagefunc{extremality_test}(copy(h))} after the
  initial setup \sage{h=\sagefunc{gmic}(2\^{ }(1/{$d$})/2)}. 
  The computation time for the construction of the field, part of the initial
  setup, does grow more rapidly with the degree, from milliseconds for small
  $d$ to seconds for $d = 60$.\par
  The number fields can be considered either abstractly or as embedded subfields 
  of an enclosing field.  When we say that the number fields are embedded into 
  the enclosing field of real numbers, this means in particular that they
  inherit the linear order from the real numbers.  To decide whether $a < b$, one computes
  sufficiently many digits of both numbers using a rigorous version of
  Newton's method; this is guaranteed to 
  terminate because the $a = b$ test can be decided by just comparing the
  coordinate vectors. 
  The program \cite{infinite-group-relaxation-code} includes a
  function \sagefunc{nice_field_values} that provides convenient access to 
  the standard facilities of SageMath that construct such an
  embedded number field.
}

We make use of the same facilities for a metaprogramming technique that
transforms the program \cite{infinite-group-relaxation-code} for testing
extremality for a function corresponding to a given parameter tuple into a
program that computes a 
description of the cell that contains the
given parameter tuple.  No code changes are necessary.

We define a class of elements\footnote{In the category-based object system 
  of SageMath, which extends Python's standard inheritance-based object system, 
  these elements are instances of the class \sagefunc{ParametricRealFieldElement}.
  Their \sage{parent}, representing the field in which the element lies, is an
  instance of the  
  class \sagefunc{ParametricRealField}.} that support the algebraic operations and
comparisons that our algorithm uses, essentially the operations of an ordered
field.  Each element stores (1) a symbolic expression\footnote{Since
  all expressions are, in fact, rational functions, we use exact
  seminumerical computations in the quotient field of a multivariate
  polynomial ring, instead of the slower and less robust general symbolic
  computation facility.} 
of the parameters in the problem, for example $x + y$ and (2) a concrete value, which
is the evaluation of this expression on the given parameter tuple, for example
$13$.  In the following, we denote elements in the form $\SymbolicElement{x+y}{13}$.
Every algebraic operation ($+$, $-$, $*$, \dots) on the elements of the class
is performed both on the symbolic expressions and on the concrete values.  
For example, if one multiplies the element $\SymbolicElement{x}{7}$ and another element
$\SymbolicElement{x+y}{13}$, one gets the element $\SymbolicElement{x^2+xy}{91}$.

When a comparison ($<$, $\leq$, $=$, \dots) takes place on elements of the
class, their concrete values are compared to compute the Boolean return value
of the comparison.  For example, the comparison
$\SymbolicElement{x^2+xy}{91} > 42$ evaluates to \emph{True}.  But we now have
a constraint on the parameters $x$ and $y$: The inequality $x^2+xy > 42$ needs to hold so
that our answer \emph{True} is correct.  We record this
constraint.\footnote{This information is recorded in the \sage{parent} of the elements.}

After a run of the algorithm, we have a description of the parameter region 
for which all the comparisons would give the same truth values as they did for the
concrete parameter tuple, and hence the algorithm would take the same
branches.  This description is typically a very long and 
highly redundant list of inequalities of rational functions in the parameters.

It is crucial to simplify the description. ``In theory'', manipulation of
inequalities describing semialgebraic sets is a solved problem of effective
real algebraic geometry.  Normal forms such as Cylindrical Algebraic
Decomposition (CAD) \cite[Chapters 5, 11]{MR2248869} are available in various
implementations, such as in the standalone QEPCAD~B
\cite{Brown:2003:QBP:968708.968710} or those integrated into CAS
software such as Maple and Mathematica, underlying these systems'
`solve' and `assume' facilities.  In computational practice, we however
observed that these systems are extremely sensitive to the number of
inequalities, rendering them unsuitable for our purposes; see \cite[section
5]{sugiyama-thesis} for a study with Maple.
We therefore roll our own implementation.  
\begin{enumerate}
\item Transform inequalities and equations of rational functions into those of
  polynomials by multiplying by denominators, and bring them in the normal
  form $p(x) < 0$ or $p(x) = 0$. 
  In the case of inequalities, this creates the extra constraint that the
  denominator takes the same sign as it does on the test point.  So this
  transformation may break cells into smaller cells.

\item Factor the polynomials $p(x)$ and record the distinct factors as equations and
  inequalities.  In the case of inequalities, this potentially breaks
  cells into smaller cells.  We can ignore the factors with even
  exponents in inequalities.

\item Reformulation--linearization:  Expand the polynomial factors in the
  standard monomial basis and replace each monomial by a new
  variable. 
  This gives a linear system of inequalities and equations and thus a
  not-necessarily-closed polyhedron in an extended variable space. 
  We use this polyhedron to represent our cell.  Indeed, its intersection with
  the algebraic variety of monomial relations is in linear bijection
  with the semialgebraic cell.

\item All of this is implemented in an incremental way.  We use the excellent
  Parma Polyhedra Library \cite{ppl-paper} via its SageMath interface written
  in Cython.  PPL is based on the double description method and supports
  not-necessarily-closed polyhedra.  It also efficiently supports adding
  inequalities dynamically and injecting a polyhedron into a higher space.
  The latter becomes necessary when a new monomial appears in some constraint.
  The PPL also has a fast implementation path for discarding redundant
  inequalities.
\end{enumerate}
In our preliminary implementation, we forgo opportunities for strengthening
this extended reformulation by McCormick inequalities, bounds propagation
etc., which would allow for further simplification.
We remark that all of these polyhedral techniques ultimately should be regarded as a
preprocessing of input for proper real-algebraic computation. They are not
strong enough on their own to provide ``minimal descriptions'' for
semialgebraic cells.  In a future version of our software, we will combine
our preprocessing technique with the CAD implementation in Mathematica.

\section{Computing the cell complex using wall-crossing BFS}
\label{s:bfs}

Define the graph of the cell complex by introducing a node for each cell and
an edge if a cell is obtained from another cell by flipping one inequality.
We compute the cell 
complex by doing a breadth-first search (BFS) in this graph.  This is a well-known
method for the case of the cells of arrangements of hyperplanes; see
\cite[chapter 24]{Goodman:1997:HDC:285869} and the references within.  The
nonlinear case poses challenges due to degeneracy and possible singularities,
which we have not completely resolved.  Our preliminary implementation uses a
heuristic numerical method to construct a point in the interior of a neighbor
cell, which will be used as the next concrete parameter tuple for re-running the
algorithm described in section~\ref{s:one-cell-metaprogramming}. 
This may fail, and so we have no guarantees that the entire parameter space is
covered by cells when the breadth-first search terminates.  This is the
weakest part of our current implementation.

\section{Automated proofs and corrections of old theorems}
\label{s:old_theorems}

\begin{figure}[h]
\centering
\begin{minipage}{.49\textwidth}
\includegraphics[width=\linewidth]{\groupcodedir/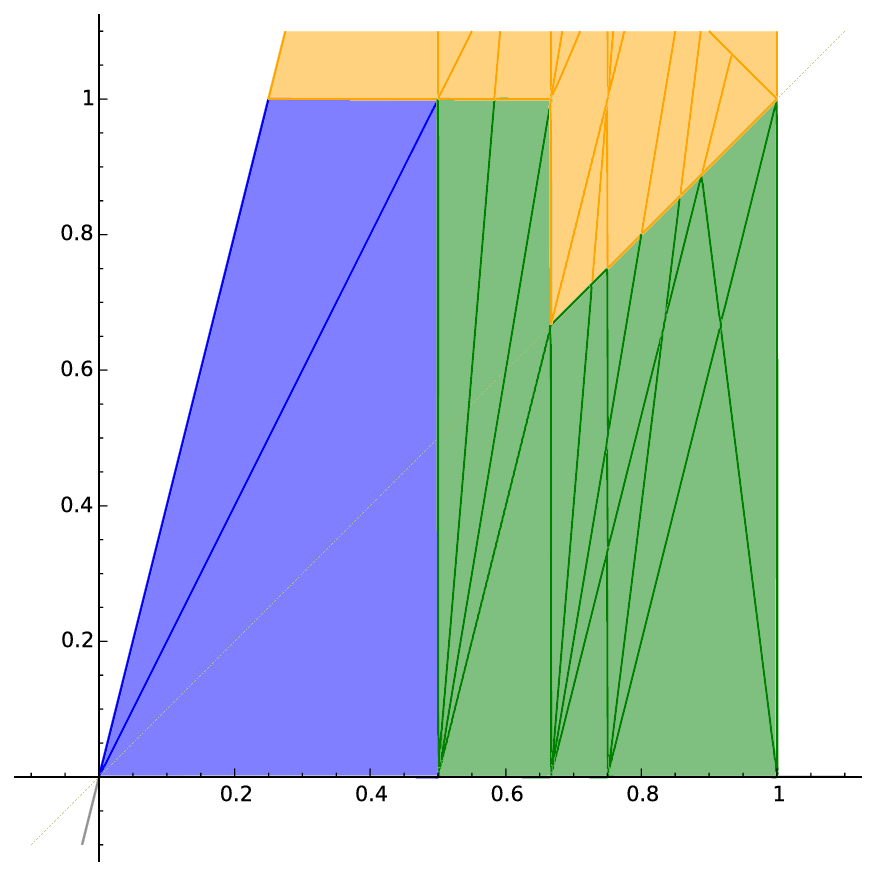}
\end{minipage}
\begin{minipage}{.49\textwidth}
\includegraphics[width=\linewidth]{\groupcodedir/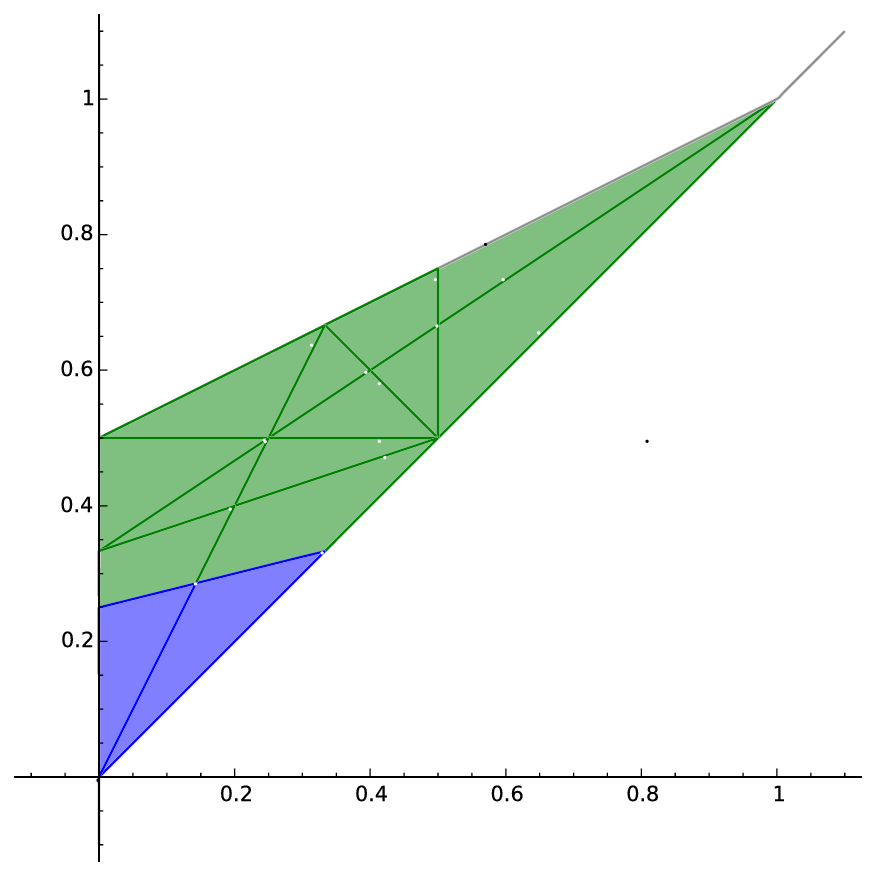}
\end{minipage}
\caption{The cell complexes of two parametric families of functions.
  \emph{Left}, \sagefunc{gj_forward_3_slope}, showing the plane of parameters
  $(\lambda_1, \lambda_2 )$ for fixed $f= 4/5$.
  \emph{Right}, \sagefunc{drlm_backward_3_slope}, showing the parameters $(f,
  \mathit{bkpt})$. 
  Cells are color-coded: `not constructible' (\emph{white}), `constructible,
  but not minimal' (\emph{yellow}), `minimal, but not
  extreme' (\emph{green}) or `extreme' (\emph{blue}).}
\label{fig:gj_drlm_complexes}
\end{figure}
Using our implementation, we verified Theorems \ref{th:gj_forward_3_slope} and
\ref{th:drlm_backward_3_slope}, as well as other theorems regarding
classical extreme functions from the literature.
Figure~\ref{fig:gj_drlm_complexes} shows 
the visualizations of the corresponding cell complexes; for additional illustrations
see Appendix~\ref{s:appendix_gj_forward_3_slope}.
Using our implementation we also investigated \autoref{th:chen_4_slope}
regarding \sagefunc{chen_4_slope} and
discovered that it is incorrect.  (See Appendix~\ref{s:appendix_chen_4_slope}
for an illustration.)
For example, 
the function with parameters $f = 7/10$, $s^+ =2$, $s^- = -4$, $\lambda_1=1/100$,
$\lambda_2 = 49/100$ satisfies the hypotheses of the theorem; however, it is not
subadditive and thus not an extreme function. On the other hand, the stated
hypotheses are also not necessary for extremality.  For example, the function
with parameters $f = 7/10$, $s^+ =2$, $s^- = -4$, $\lambda_1=1/10$, $\lambda_2 = 1/10$
does not satisfy the hypotheses, however it is extreme. 
We omit a statement of corrected hypotheses that we found using our code.

We also investigated another family of functions, the so-called
$\mathrm{CPL}_3^=$ functions, introduced by the systematic
study by Miller, Li, and Richard \cite{Miller-Li-Richard2008}.  Their method
can be regarded as a predecessor of our method, albeit one that led to an
error-prone manual case analysis (and human-generated proofs). 
Though our general method can be applied directly
, we developed a specialized version of our code that follows Miller, Li, and
Richard's method to allow a direct comparison. 
This revealed 
mistakes in \cite{Miller-Li-Richard2008}, as we report
in Appendix~\ref{s:appendix_cpl3}.

\section{Computer-assisted discovery of new theorems}
\label{s:new-theorems}

In \cite{koeppe-zhou:extreme-search}, the authors conducted a systematic
computer-based search for extreme functions on the grids $\frac1q\Z$ for
values of $q$ up to 30.  This resulted in a large catalog of extreme functions
that are ``sporadic'' in the sense that they do not belong to any parametric
family described in the literature.

Our goal is to automatically embed these functions into parametric families
and to automatically prove theorems about their extremality.  In this section,
we report on cases that have been done successfully with our preliminary
implementation; the process is not completely automatic yet.

\begin{figure}[t]
\centering
\includegraphics[width=0.45\textwidth]{\groupcodedir/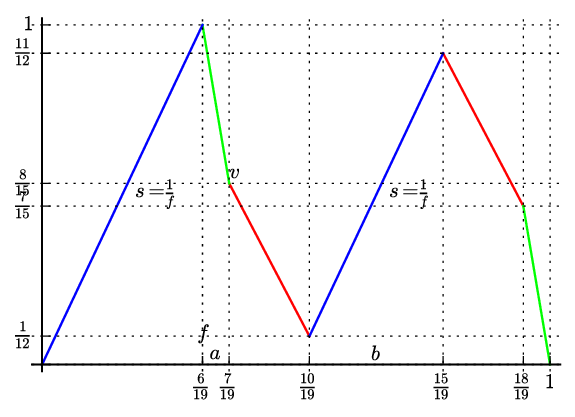}
\quad
\includegraphics[width=0.45\textwidth]{\groupcodedir/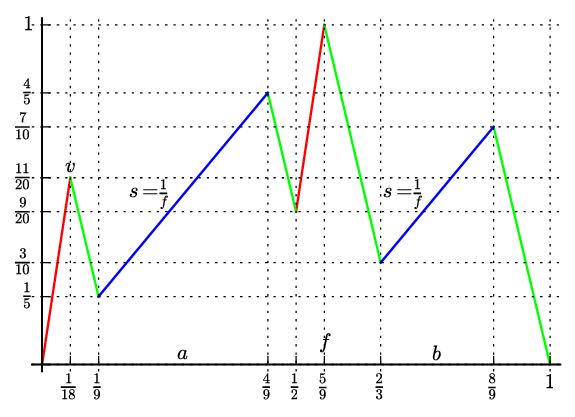}
\caption{Two new parametric families of extreme functions. \emph{Left},
  \sagefunc{kzh_3_slope_param_extreme_1}. \emph{Right},
  \sagefunc{kzh_3_slope_param_extreme_2}.}
\label{fig:plot_kzh_3_slope_param_extreme_1}
\label{fig:plot_kzh_3_slope_param_extreme_2}
\end{figure}

We picked an interesting-looking $3$-slope extreme function found by our
computer-based search on the grid $\frac{1}{q}\Z$. 
We then introduced parameters $f, a, b,v$ to describe a preliminary parametric
family that we denote by \sage{param\_3\_slope\_1}.
In the concrete function that we started from, these parameters take the
values $\frac{6}{19},\frac{1}{19},\frac{5}{19},\frac{8}{15}$; 
see
Figure~\ref{fig:plot_kzh_3_slope_param_extreme_1} (left).  So  $a$ denotes the
length of the first interval right to $f$, $b$ denotes the length of interval
centered at $(1+f)/2$ and $v = \pi(f+a)$.  By this choice of parameters, the
function automatically satisfies the equations corresponding to the symmetry
conditions. 
Next we run the parametric version of the minimality test algorithm.  It
computes \emph{True}, and as a side-effect computes a description of the cell
in which the minimality test is the same.
{\small
\begin{verbatim}
sage: K.<f,a,b,v>=ParametricRealField([6/19,1/19,5/19,8/15])
sage: h = param_3_slope_1(f,a,b,v)
sage: minimality_test(h)
True
sage: K._eq_factor
{-f^2*v + 3*f*b*v + f^2 + f*a - 3*f*b - 3*a*b - f*v + b}
\end{verbatim}}
In particular, the above line shows that it has discovered one nonlinear equation that holds in the cell
corresponding to the minimality proof of the concrete function.  
We use this equation, quadratic in $f$ and multilinear in the other
parameters, to eliminate one of the parameters.\footnote{We plan to automate this in
a future version of our software.}
This gives our parametric family \sagefunc{kzh_3_slope_param_extreme_1},
which depends only on parameters $f, a, b$.
See Appendix~\ref{s:appendix_kzh_3_slope_param_definitions} to see how this family
is defined in the code; the definition of the parametric family is the only
input to our algorithm. 
Cells with respect to this family
will be full-dimensional; this helps to satisfy a current implementation
restriction of our software.
Indeed, re-running the algorithm yields the following simplified description
of the cell in which the concrete parameter tuple lies.
{\small
\begin{verbatim}
  3*f + 4*a - b - 1 < 0                    -a < 0           
  -f^2 - f*a + 3*f*b + 3*a*b - b < 0       -f + b < 0       
  f*a - 3*a*b - f + b < 0                  -f - 3*b + 1 < 0 
  -f^2*a + 3*f*a*b - 3*a*b - f + b < 0
\end{verbatim}}
We then compute the cell complex by BFS as described in section~\ref{s:bfs};
see Appendix~\ref{s:appendix_kzh_3_slope_param_extreme_1} for an illustration.
By inspection, we observe that the collection of the cells for which the
function is extreme happens to be a convex polytope (this is not
guaranteed). We discard the inequalities that appear twice and thus
describe inner walls of the complex.  By inspection, we discard 
nonlinear inequalities that are redundant.  We obtain a description of the
union of the cells for which the function is extreme as a convex polytope.
We obtain the following:
\begin{theorem} 
Let $f \in (0,1)$ and $a, b \in \R$ such that \[0 \leq a, \quad 0 \leq b \leq
  f \text{ and } 3f+4a -b -1 \leq0.\]
The periodic, piecewise linear \sagefunc{kzh_3_slope_param_extreme_1} function $\pi\colon \R/\Z \to \R$ defined as follows is extreme. The function $\pi$ has breakpoints at \[0, f, f+a, \frac{1+f-b}{2}, \frac{1+f+b}{2}, 1-a, 1.\] The values at breakpoints are given by $\pi(0) = \pi(1) =0$, $\pi(f+a)=1-\pi(1-a)=v$ and $\pi(\frac{1+f-b}{2})=1-\pi(\frac{1+f+b}{2})=\frac{f-b}{2f}$,  where  $v= \frac{f^2+fa-3fb-3ab+b}{f^2+f-3bf}$.
\end{theorem}

A similar process leads to a theorem about the family
\sagefunc{kzh_3_slope_param_extreme_2} shown in
Figure~\ref{fig:plot_kzh_3_slope_param_extreme_2} (right). We omit the
statement of the theorem.
Appendix~\ref{s:appendix_kzh_3_slope_param_extreme_2} shows a visualization of
the cell complex.

{\scriptsize
\providecommand\ISBN{ISBN }
\bibliographystyle{../../amsabbrvurl}
\bibliography{../../bib/MLFCB_bib}
}
\clearpage
\appendix
\section{Appendix}

\subsection{Additional illustrations for \sage{gj\_forward\_3\_slope}}
\label{s:appendix_gj_forward_3_slope}

In Figure~\ref{fig:gj_drlm_complexes} in section~\ref{s:old_theorems} we
showed the slice of the cell complex for fixed parameter $f=4/5$.
In Figure~\ref{fig:gj_forward_3_slope_f_lambda_1} below, we show two other views
on the parameter space. 

\begin{figure}[h]
\centering
\begin{minipage}{.5\textwidth}
\includegraphics[width=\linewidth]{\groupcodedir/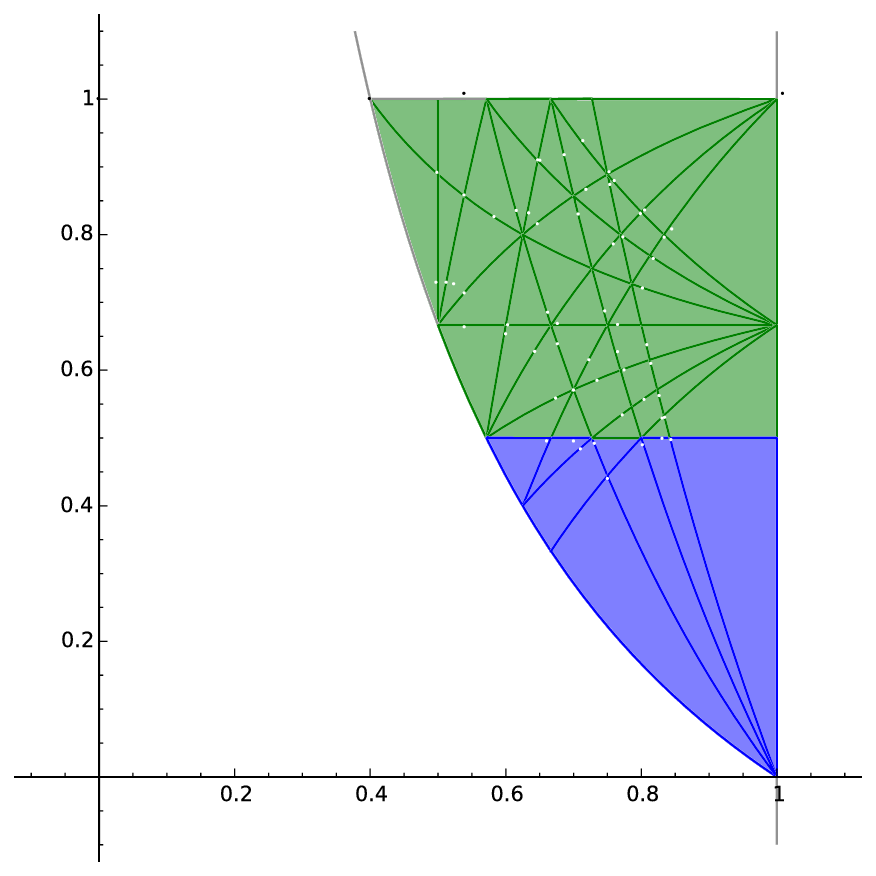}
\end{minipage}%
\begin{minipage}{.5\textwidth}
\includegraphics[width=\linewidth]{\groupcodedir/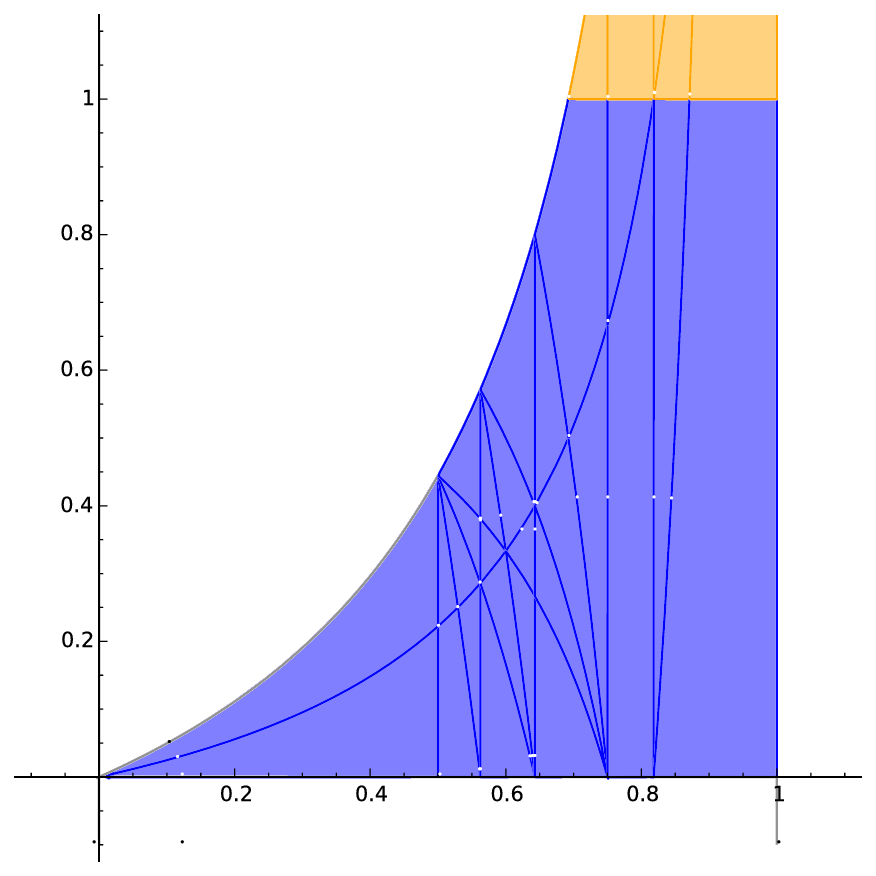}
\end{minipage}
\caption{The cell complex of the parametric family
  \sagefunc{gj_forward_3_slope}.
  \emph{Left}, showing the plane of parameters $(f, \lambda_1 )$ for fixed $\lambda_2 = 2/3$; 
  \emph{Right}, showing the plane of parameters $(f, \lambda_2 )$ for fixed $\lambda_1 = 4/9$. 
  Cells are color-coded: `not constructible' (\emph{white}), `constructible,
  but not minimal' (\emph{yellow}), `minimal, but not
  extreme' (\emph{green}) or `extreme' (\emph{blue}).}
\label{fig:gj_forward_3_slope_f_lambda_1}
\end{figure}

\clearpage
\subsection{Cell complex for \sage{chen\_4\_slope}}
\label{s:appendix_chen_4_slope}

\begin{figure}[h]
\centering
\begin{minipage}{.49\textwidth}
\includegraphics[width=\linewidth]{\groupcodedir/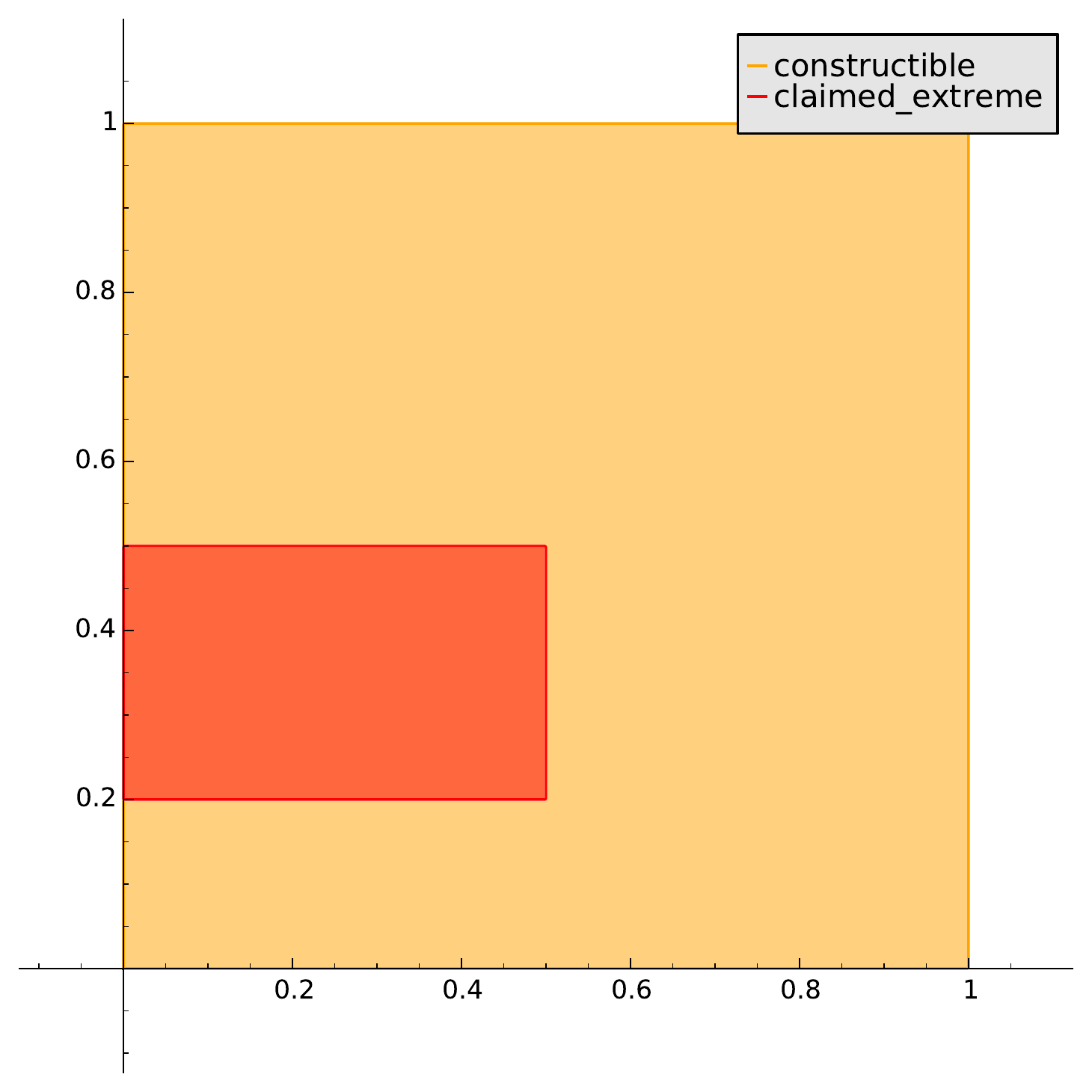}
\end{minipage}
\begin{minipage}{.49\textwidth}
\includegraphics[width=\linewidth]{\groupcodedir/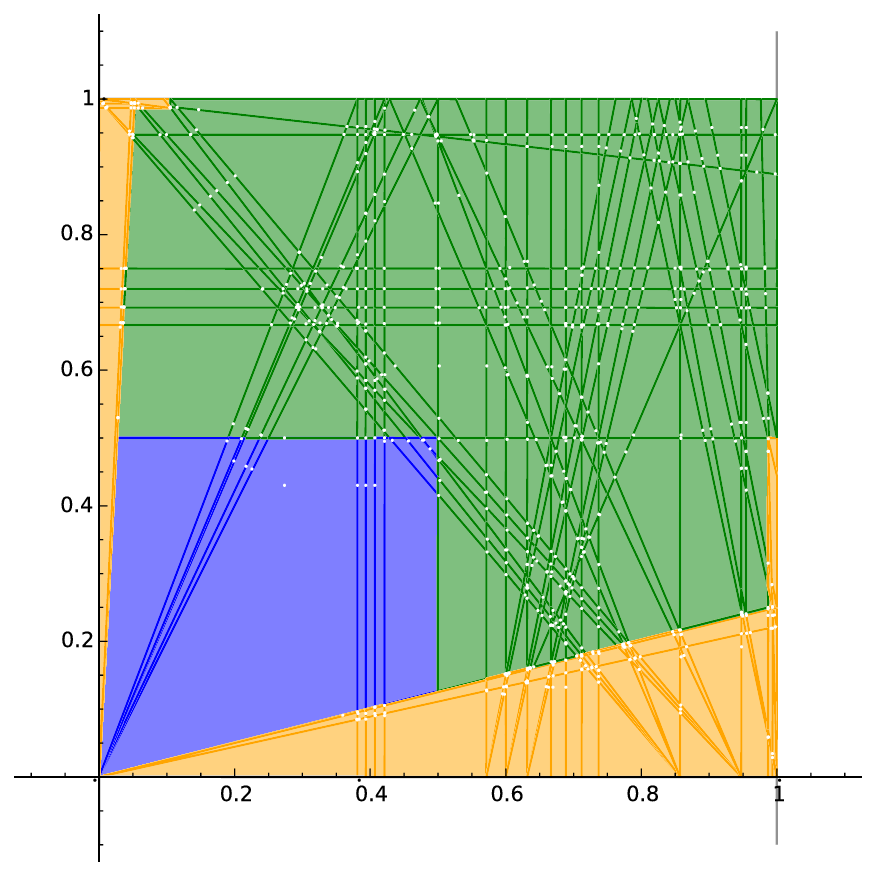}
\end{minipage}
\caption{The extreme region of \sagefunc{chen_4_slope} claimed in the
  literature is incorrect. Parameters $(\lambda_1, \lambda_2)$; fixed $f=7/10$,
  $s^+=2$, $s^-=-4$.
  \emph{Left}, the incorrect hypotheses of Theorem~\ref{th:chen_4_slope}
  (\emph{orange}) within the region of constructibility of the function (\emph{yellow}).
  \emph{Right}, the cell complex computed by our implementation.  
  Cells are color-coded: `constructible,
  but not minimal' (\emph{yellow}), `minimal, but not
  extreme' (\emph{green}) or `extreme' (\emph{blue}).}
\label{fig:chen_4_slope_complex}
\end{figure}

Figure~\ref{fig:chen_4_slope_complex} illustrates our remarks in
section~\ref{s:old_theorems} regarding the extremality conditions in Chen's
theorem (Theorem~\ref{th:chen_4_slope}) about his \sagefunc{chen_4_slope}
family.  The parameter space is 5-dimensional; the figures show a
2-dimensional slice corresponding to varying parameters $\lambda_1$ and
$\lambda_2$ and fixed parameters $f$, $s^+$, $s^-$. 

\clearpage
\subsection{Miller, Li, and Richard's $\text{CPL}_3^=$-extreme functions}
\label{s:appendix_cpl3}

Miller, Li and Richard \cite{Miller-Li-Richard2008} derived extreme functions
using an approximate lifting scheme proposed in \cite{Richard-Li-Miller-2009:Approximate-Liftings}. 
The authors studied a family of continuous piecewise linear functions $\phi$, called $\text{CPL}_3^=$-lifting functions, that have $4$ parameters $r_0, z_1, \theta_1, \theta_2$. ($r_0$ is our $f$.)
Given $(r_0, z_1) \in \R^2$  such that $r_0 +4z_1 \leq 1$, the parameters $(\theta_1, \theta_2)$ that define a superadditive $\text{CPL}_3^=$-lifting function $\phi$ belong to a certain polytope $P\Theta_3^=$.
They investigated the $\text{CPL}_3^=$-group functions $\pi$ obtained by converting the $\text{CPL}_3^=$-lifting functions $\phi$ that correspond to extreme points of the polytope $P\Theta_3^=$. These functions $\pi$ are potentially strong valid inequalities for the infinite group problem, but they are not always extreme. As a result of a manual inspection which required extensive case analysis, the authors summarized in  \cite[Table 5]{Miller-Li-Richard2008} the conditions on the parameters $r_0, z_1, \theta_1, \theta_2$ for which the corresponding $\text{CPL}_3^=$-group functions $\pi$  are extreme for the infinite group problem. 
These extreme functions $\pi$ belong to $14$ parametric families of two- and three-slope continuous piecewise linear functions, among which some had not been discovered before. 
The objective of this section is to verify and reproduce the results stated in \cite{Miller-Li-Richard2008} by an automated discovery process.

In order to get the description of the polytope $P\Theta_3^=$, one writes out all the superadditivity constraints $\phi(x)+\phi(y)\leq\phi(z)$ where $z=x+y$, $(x, y)$ or $(x,z)$ is a pair of breakpoints of the $\text{CPL}_3^=$-lifting function $\phi$. (Or equivalently all the subadditivity constraints on the breakpoints of the $\text{CPL}_3^=$-group functions $\pi$.)
For example, at the breakpoint $r_0+z_1$, $\phi$ takes value $\theta_1$, so the superadditivity constraint $\phi(r_0+z_1)+ \phi(r_0+z_1) \leq \phi(2r_0+2z_1)$ reduces to $2\theta_1 \leq \phi(2r_0+2z_1)$. This inequality is satisfied by any $(\theta_1, \theta_2) \in P\Theta_3^=$. To determine the value $\phi(2r_0+2z_1)$ in terms of $\theta_1$ and $\theta_2$, one needs to know which linear piece of $\phi$ the value $2r_0+2z_1$ falls into. 
So the form of an inequality that defines $P\Theta_3^=$ depends on how the breakpoints of $\phi$ are arranged. We would first compute the cells of $(r_0, z_1)$ where the arrangement of breakpoints stays the same. 
That is, if for $\phi$ corresponding to a given $(r_0^*, z_1^*)$, the sum (or difference) of its $i$-th and $j$-th breakpoints is contained in its $k$-th linear piece, then same holds for any $(r_0, z_1)$ inside the cell which contains $(r_0^*, z_1^*)$.
In consequence, the inequalities defining $P\Theta_3^=$ take a fixed form for any $(r_0, z_1)$ inside a cell, but differ from cell to cell.  Figure~\ref{fig:region_87} illustrates these cells that are obtained by running the following code:
\begin{verbatim}
sage: regions = regions_r0_z1_from_arrangement_of_bkpts()
sage: g = plot_cpl_components(regions)
sage: g.show(xmin=0, xmax=1, ymin=0, ymax=1/4)
\end{verbatim}
We observe that the result agrees with \cite[Table 1]{Miller-Li-Richard2008}. 

We then study these cells one by one. Assume that $(r_0, z_1)$ is restricted to a cell where the arrangement of breakpoints is combinatorially the same.
The superadditivity constraints on pairs of breakpoints of $\phi$ give a set of $N=52$ linear inequalities that define $P\Theta_3^=$. This description includes many redundant inequalities. To save time on the case analysis later, a simplified description of $P\Theta_3^=$ which consists of $N=9$ linear inequalities 
is obtained in \cite[section 3.1]{Miller-Li-Richard2008}. This step is not necessary in our automated study, because the code can afford $N=52$.

The value of $(\theta_1, \theta_2)$  in $(r_0, z_1, \theta_1, \theta_2)$ that
gives an extreme function $\pi$ must be an extreme point $(\theta_1,
\theta_2)$ of the polytope $P\Theta_3^=$. To obtain the extreme point
$(\theta_1, \theta_2)$, we pick $2$  inequalities to be tight in the
description of $P\Theta_3^=$. Suppose that the $i$-th and the $j$-th
inequalities of $P\Theta_3^=$ are tight, for $1 \leq i < j \leq N$. One can
solve for $(\theta_1, \theta_2)$ symbolically in terms of $(r_0, z_1)$. Then,
we plug in this solution $(\theta_1, \theta_2)$ in $\pi$ to get a
$\text{CPL}_3^=$-group function $\pi$ with parameters $r_0$ and $z_1$.  Run
the BFS (always restricted to a  $(r_0, z_1)$ cell with same arrangement of
breakpoints) to find regions where $\pi$ is extreme (\emph{blue}); is minimal but not
extreme (\emph{green}); is not minimal (\emph{yellow}) (not minimal implies some
subadditivity constraint of $\pi$ is violated, so the solution $(\theta_1,
\theta_2)$ does not satisfy all other inequalities of $P\Theta_3^=$, and thus
$(\theta_1, \theta_2)$  is not an extreme point of the polytope at this $(r_0,
z_1)$.) Repeat this process for all $1\leq i< j\leq N$, and for all
arrangement cells. For each different expression of $(\theta_1, \theta_2)$
that has an extreme (\emph{blue}) region, plot the results on a diagram. See  \autoref{fig:cpl_theta_i}. These diagrams are obtained by running the following code:
\begin{verbatim}
# 87 regions in Figure 8. 
# For each region, find theta solutions, then subdivide into
# smaller components by running bfs with parametric field.
sage: regions = cpl_regions_with_thetas_and_components()
# Gather the blue components that correspond to the same
# expression of theta together. 
sage: thetas_and_regions = cpl_thetas_and_regions_extreme(regions)
# Get diagrams "cpl_ext_theta_i" that show only blue regions.
sage: save_cpl_extreme_theta_regions(thetas_and_regions)
# Get diagrams "cpl_theta_i", that show a bit more.
sage: thetas_and_components = {}
sage: for theta in thetas_and_regions.keys():
....:     components = cpl_regions_fix_theta(regions, theta)
....:     thetas_and_components[theta]=components
sage: save_cpl_extreme_theta_regions(thetas_and_components)
\end{verbatim}
We compare the diagrams (and their inequality descriptions) with the
conditions for extremality stated in \cite[Table 5]{Miller-Li-Richard2008}. We
verify that the code is successful in finding all extreme regions claimed in
the literature.\footnote{Again
we note that some heuristics and manual inspection were necessary, in particular
to deal with some redundant nonlinear inequalities, which our current
implementation is not able to remove; see the discussion at the end of
\autoref{s:semialgebraic-cells}.}  
Furthermore, we observe that for the case `k' the condition for extremality $r_0 \leq 2z_1$
and $r_0 + 5z_1=1$ stated in \cite[Table 5]{Miller-Li-Richard2008} is
incorrect; it should be $r_0 \geq z_1$ and $r_0 + 5z_1=1$. 
We remark that in the course of adding these functions
(\sage{mlr\_cpl3\_}\dots) to the Electronic
Compendium~\cite{infinite-group-relaxation-code}, Sugiyama
\cite{sugiyama-thesis} had already found another, unrelated mistake: \cite[Table
3]{Miller-Li-Richard2008} shows incorrect slope values $s_3$ of case `l' and $s_4$
of case `p'.

\input{cpl_extreme_cases_diagrams.tex}

\clearpage
\subsection{Definition of the new functions}
\label{s:appendix_kzh_3_slope_param_definitions}

The following Sage code defines the new families of functions. 

{\scriptsize
\begin{verbatim}
def kzh_3_slope_param_extreme_1(f=6/19, a=1/19, b=5/19, field=None, conditioncheck=True):
    if not bool(0 < f < f+a < (1+f-b)/2 < (1+f+b)/2 < 1-a < 1):
        raise ValueError, "Bad parameters. Unable to construct the function."
    if conditioncheck:
        if not bool(0 <= a and  0 <= b <= f and 3*f+4*a-b-1 <= 0):
            logging.info("Conditions for extremality are NOT satisfied.")
        else:
            logging.info("Conditions for extremality are satisfied.")
    v = (f*f+f*a-3*f*b-3*a*b+b)/(f*f+f-3*f*b)
    bkpts = [0, f, f+a, (1+f-b)/2, (1+f+b)/2, 1-a, 1]
    values = [0, 1, v, (f-b)/2/f, (f+b)/2/f, 1-v, 0]
    return piecewise_function_from_breakpoints_and_values(bkpts, values, field=field)

def kzh_3_slope_param_extreme_2(f=5/9, a=3/9, b=2/9, field=None, conditioncheck=True):
    if not bool(0 < a < f < 1 and 0 < b < f):
        raise ValueError, "Bad parameters. Unable to construct the function."
    if conditioncheck:
        if not bool(b <= a and f <= a + b and f <= (1+a-b)/2):
            logging.info("Conditions for extremality are NOT satisfied.")
        else:
            logging.info("Conditions for extremality are satisfied.")
    v = (f*(f-a+b-2)-a*b+2*a)/(f+b-1)/f/4;
    bkpts = [0, (f-a)/4, (f-a)/2, (f+a)/2, f-(f-a)/4, f, \
             (1+f-b)/2, (1+f+b)/2, 1]
    values = [0, v, (f-a)/f/2, (f+a)/f/2, 1-v, 1, (f-b)/f/2, (f+b)/f/2, 0]
    return piecewise_function_from_breakpoints_and_values(bkpts, values, field=field)

\end{verbatim}
}

\clearpage
\subsection{Cell complex for \sage{kzh\_3\_slope\_param\_extreme\_1}}
\label{s:appendix_kzh_3_slope_param_extreme_1}

Figure~\ref{fig:kzh_3_slope_param_extreme_1_slices} illustrates the
3-dimensional cell complex for the new family
\sagefunc{kzh_3_slope_param_extreme_1}, introduced in
section~\ref{s:new-theorems}, by showing slices for various fixed values of
the parameter~$f$.  

\input{param_3s_ext_1_f_diagrams.tex}

\clearpage
\subsection{Cell complex for \sage{kzh\_3\_slope\_param\_extreme\_2}}
\label{s:appendix_kzh_3_slope_param_extreme_2}

Figure~\ref{fig:kzh_3_slope_param_extreme_1_slices} illustrates the
3-dimensional cell complex for the new family
\sagefunc{kzh_3_slope_param_extreme_2}, introduced in
section~\ref{s:new-theorems}, by showing slices for various fixed values of
the parameter~$f$.  

\input{param_3s_ext_2_f_diagrams.tex}

\end{document}

%% file: sage-commands.tex
\pgfkeyssetvalue{/sagefunc/bccz_counterexample}{\href{\githubsearchurl?q=\%22def+bccz_counterexample(\%22}{\sage{bccz\underscore{}counterexample}}}
\pgfkeyssetvalue{/sagefunc/bcdsp_arbitrary_slope}{\href{\githubsearchurl?q=\%22def+bcdsp_arbitrary_slope(\%22}{\sage{bcdsp\underscore{}arbitrary\underscore{}slope}}}
\pgfkeyssetvalue{/sagefunc/bhk_gmi_irrational}{\href{\githubsearchurl?q=\%22def+bhk_gmi_irrational(\%22}{\sage{bhk\underscore{}gmi\underscore{}irrational}}}
\pgfkeyssetvalue{/sagefunc/bhk_irrational}{\href{\githubsearchurl?q=\%22def+bhk_irrational(\%22}{\sage{bhk\underscore{}irrational}}}
\pgfkeyssetvalue{/sagefunc/bhk_slant_irrational}{\href{\githubsearchurl?q=\%22def+bhk_slant_irrational(\%22}{\sage{bhk\underscore{}slant\underscore{}irrational}}}
\pgfkeyssetvalue{/sagefunc/chen_4_slope}{\href{\githubsearchurl?q=\%22def+chen_4_slope(\%22}{\sage{chen\underscore{}4\underscore{}slope}}}
\pgfkeyssetvalue{/sagefunc/dg_2_step_mir}{\href{\githubsearchurl?q=\%22def+dg_2_step_mir(\%22}{\sage{dg\underscore{}2\underscore{}step\underscore{}mir}}}
\pgfkeyssetvalue{/sagefunc/dg_2_step_mir_limit}{\href{\githubsearchurl?q=\%22def+dg_2_step_mir_limit(\%22}{\sage{dg\underscore{}2\underscore{}step\underscore{}mir\underscore{}limit}}}
\pgfkeyssetvalue{/sagefunc/drlm_2_slope_limit}{\href{\githubsearchurl?q=\%22def+drlm_2_slope_limit(\%22}{\sage{drlm\underscore{}2\underscore{}slope\underscore{}limit}}}
\pgfkeyssetvalue{/sagefunc/drlm_3_slope_limit}{\href{\githubsearchurl?q=\%22def+drlm_3_slope_limit(\%22}{\sage{drlm\underscore{}3\underscore{}slope\underscore{}limit}}}
\pgfkeyssetvalue{/sagefunc/drlm_backward_3_slope}{\href{\githubsearchurl?q=\%22def+drlm_backward_3_slope(\%22}{\sage{drlm\underscore{}backward\underscore{}3\underscore{}slope}}}
\pgfkeyssetvalue{/sagefunc/gj_2_slope}{\href{\githubsearchurl?q=\%22def+gj_2_slope(\%22}{\sage{gj\underscore{}2\underscore{}slope}}}
\pgfkeyssetvalue{/sagefunc/gj_2_slope_repeat}{\href{\githubsearchurl?q=\%22def+gj_2_slope_repeat(\%22}{\sage{gj\underscore{}2\underscore{}slope\underscore{}repeat}}}
\pgfkeyssetvalue{/sagefunc/gj_forward_3_slope}{\href{\githubsearchurl?q=\%22def+gj_forward_3_slope(\%22}{\sage{gj\underscore{}forward\underscore{}3\underscore{}slope}}}
\pgfkeyssetvalue{/sagefunc/gmic}{\href{\githubsearchurl?q=\%22def+gmic(\%22}{\sage{gmic}}}
\pgfkeyssetvalue{/sagefunc/hildebrand_2_sided_discont_1_slope_1}{\href{\githubsearchurl?q=\%22def+hildebrand_2_sided_discont_1_slope_1(\%22}{\sage{hildebrand\underscore{}2\underscore{}sided\underscore{}discont\underscore{}1\underscore{}slope\underscore{}1}}}
\pgfkeyssetvalue{/sagefunc/hildebrand_2_sided_discont_2_slope_1}{\href{\githubsearchurl?q=\%22def+hildebrand_2_sided_discont_2_slope_1(\%22}{\sage{hildebrand\underscore{}2\underscore{}sided\underscore{}discont\underscore{}2\underscore{}slope\underscore{}1}}}
\pgfkeyssetvalue{/sagefunc/hildebrand_5_slope_22_1}{\href{\githubsearchurl?q=\%22def+hildebrand_5_slope_22_1(\%22}{\sage{hildebrand\underscore{}5\underscore{}slope\underscore{}22\underscore{}1}}}
\pgfkeyssetvalue{/sagefunc/hildebrand_5_slope_24_1}{\href{\githubsearchurl?q=\%22def+hildebrand_5_slope_24_1(\%22}{\sage{hildebrand\underscore{}5\underscore{}slope\underscore{}24\underscore{}1}}}
\pgfkeyssetvalue{/sagefunc/hildebrand_5_slope_28_1}{\href{\githubsearchurl?q=\%22def+hildebrand_5_slope_28_1(\%22}{\sage{hildebrand\underscore{}5\underscore{}slope\underscore{}28\underscore{}1}}}
\pgfkeyssetvalue{/sagefunc/hildebrand_discont_3_slope_1}{\href{\githubsearchurl?q=\%22def+hildebrand_discont_3_slope_1(\%22}{\sage{hildebrand\underscore{}discont\underscore{}3\underscore{}slope\underscore{}1}}}
\pgfkeyssetvalue{/sagefunc/kf_n_step_mir}{\href{\githubsearchurl?q=\%22def+kf_n_step_mir(\%22}{\sage{kf\underscore{}n\underscore{}step\underscore{}mir}}}
\pgfkeyssetvalue{/sagefunc/kzh_10_slope_1}{\href{\githubsearchurl?q=\%22def+kzh_10_slope_1(\%22}{\sage{kzh\underscore{}10\underscore{}slope\underscore{}1}}}
\pgfkeyssetvalue{/sagefunc/kzh_28_slope_1}{\href{\githubsearchurl?q=\%22def+kzh_28_slope_1(\%22}{\sage{kzh\underscore{}28\underscore{}slope\underscore{}1}}}
\pgfkeyssetvalue{/sagefunc/kzh_28_slope_2}{\href{\githubsearchurl?q=\%22def+kzh_28_slope_2(\%22}{\sage{kzh\underscore{}28\underscore{}slope\underscore{}2}}}
\pgfkeyssetvalue{/sagefunc/kzh_3_slope_param_extreme_1}{\href{\githubsearchurl?q=\%22def+kzh_3_slope_param_extreme_1(\%22}{\sage{kzh\underscore{}3\underscore{}slope\underscore{}param\underscore{}extreme\underscore{}1}}}
\pgfkeyssetvalue{/sagefunc/kzh_3_slope_param_extreme_2}{\href{\githubsearchurl?q=\%22def+kzh_3_slope_param_extreme_2(\%22}{\sage{kzh\underscore{}3\underscore{}slope\underscore{}param\underscore{}extreme\underscore{}2}}}
\pgfkeyssetvalue{/sagefunc/kzh_4_slope_param_extreme_1}{\href{\githubsearchurl?q=\%22def+kzh_4_slope_param_extreme_1(\%22}{\sage{kzh\underscore{}4\underscore{}slope\underscore{}param\underscore{}extreme\underscore{}1}}}
\pgfkeyssetvalue{/sagefunc/kzh_5_slope_fulldim_1}{\href{\githubsearchurl?q=\%22def+kzh_5_slope_fulldim_1(\%22}{\sage{kzh\underscore{}5\underscore{}slope\underscore{}fulldim\underscore{}1}}}
\pgfkeyssetvalue{/sagefunc/kzh_5_slope_fulldim_2}{\href{\githubsearchurl?q=\%22def+kzh_5_slope_fulldim_2(\%22}{\sage{kzh\underscore{}5\underscore{}slope\underscore{}fulldim\underscore{}2}}}
\pgfkeyssetvalue{/sagefunc/kzh_5_slope_fulldim_3}{\href{\githubsearchurl?q=\%22def+kzh_5_slope_fulldim_3(\%22}{\sage{kzh\underscore{}5\underscore{}slope\underscore{}fulldim\underscore{}3}}}
\pgfkeyssetvalue{/sagefunc/kzh_5_slope_fulldim_4}{\href{\githubsearchurl?q=\%22def+kzh_5_slope_fulldim_4(\%22}{\sage{kzh\underscore{}5\underscore{}slope\underscore{}fulldim\underscore{}4}}}
\pgfkeyssetvalue{/sagefunc/kzh_5_slope_fulldim_5}{\href{\githubsearchurl?q=\%22def+kzh_5_slope_fulldim_5(\%22}{\sage{kzh\underscore{}5\underscore{}slope\underscore{}fulldim\underscore{}5}}}
\pgfkeyssetvalue{/sagefunc/kzh_5_slope_fulldim_covers_1}{\href{\githubsearchurl?q=\%22def+kzh_5_slope_fulldim_covers_1(\%22}{\sage{kzh\underscore{}5\underscore{}slope\underscore{}fulldim\underscore{}covers\underscore{}1}}}
\pgfkeyssetvalue{/sagefunc/kzh_5_slope_fulldim_covers_2}{\href{\githubsearchurl?q=\%22def+kzh_5_slope_fulldim_covers_2(\%22}{\sage{kzh\underscore{}5\underscore{}slope\underscore{}fulldim\underscore{}covers\underscore{}2}}}
\pgfkeyssetvalue{/sagefunc/kzh_5_slope_fulldim_covers_3}{\href{\githubsearchurl?q=\%22def+kzh_5_slope_fulldim_covers_3(\%22}{\sage{kzh\underscore{}5\underscore{}slope\underscore{}fulldim\underscore{}covers\underscore{}3}}}
\pgfkeyssetvalue{/sagefunc/kzh_5_slope_fulldim_covers_4}{\href{\githubsearchurl?q=\%22def+kzh_5_slope_fulldim_covers_4(\%22}{\sage{kzh\underscore{}5\underscore{}slope\underscore{}fulldim\underscore{}covers\underscore{}4}}}
\pgfkeyssetvalue{/sagefunc/kzh_5_slope_fulldim_covers_5}{\href{\githubsearchurl?q=\%22def+kzh_5_slope_fulldim_covers_5(\%22}{\sage{kzh\underscore{}5\underscore{}slope\underscore{}fulldim\underscore{}covers\underscore{}5}}}
\pgfkeyssetvalue{/sagefunc/kzh_5_slope_fulldim_covers_6}{\href{\githubsearchurl?q=\%22def+kzh_5_slope_fulldim_covers_6(\%22}{\sage{kzh\underscore{}5\underscore{}slope\underscore{}fulldim\underscore{}covers\underscore{}6}}}
\pgfkeyssetvalue{/sagefunc/kzh_5_slope_q22_f10_1}{\href{\githubsearchurl?q=\%22def+kzh_5_slope_q22_f10_1(\%22}{\sage{kzh\underscore{}5\underscore{}slope\underscore{}q22\underscore{}f10\underscore{}1}}}
\pgfkeyssetvalue{/sagefunc/kzh_5_slope_q22_f10_2}{\href{\githubsearchurl?q=\%22def+kzh_5_slope_q22_f10_2(\%22}{\sage{kzh\underscore{}5\underscore{}slope\underscore{}q22\underscore{}f10\underscore{}2}}}
\pgfkeyssetvalue{/sagefunc/kzh_5_slope_q22_f10_3}{\href{\githubsearchurl?q=\%22def+kzh_5_slope_q22_f10_3(\%22}{\sage{kzh\underscore{}5\underscore{}slope\underscore{}q22\underscore{}f10\underscore{}3}}}
\pgfkeyssetvalue{/sagefunc/kzh_5_slope_q22_f10_4}{\href{\githubsearchurl?q=\%22def+kzh_5_slope_q22_f10_4(\%22}{\sage{kzh\underscore{}5\underscore{}slope\underscore{}q22\underscore{}f10\underscore{}4}}}
\pgfkeyssetvalue{/sagefunc/kzh_5_slope_q22_f2_1}{\href{\githubsearchurl?q=\%22def+kzh_5_slope_q22_f2_1(\%22}{\sage{kzh\underscore{}5\underscore{}slope\underscore{}q22\underscore{}f2\underscore{}1}}}
\pgfkeyssetvalue{/sagefunc/kzh_6_slope_1}{\href{\githubsearchurl?q=\%22def+kzh_6_slope_1(\%22}{\sage{kzh\underscore{}6\underscore{}slope\underscore{}1}}}
\pgfkeyssetvalue{/sagefunc/kzh_6_slope_fulldim_covers_1}{\href{\githubsearchurl?q=\%22def+kzh_6_slope_fulldim_covers_1(\%22}{\sage{kzh\underscore{}6\underscore{}slope\underscore{}fulldim\underscore{}covers\underscore{}1}}}
\pgfkeyssetvalue{/sagefunc/kzh_6_slope_fulldim_covers_2}{\href{\githubsearchurl?q=\%22def+kzh_6_slope_fulldim_covers_2(\%22}{\sage{kzh\underscore{}6\underscore{}slope\underscore{}fulldim\underscore{}covers\underscore{}2}}}
\pgfkeyssetvalue{/sagefunc/kzh_6_slope_fulldim_covers_3}{\href{\githubsearchurl?q=\%22def+kzh_6_slope_fulldim_covers_3(\%22}{\sage{kzh\underscore{}6\underscore{}slope\underscore{}fulldim\underscore{}covers\underscore{}3}}}
\pgfkeyssetvalue{/sagefunc/kzh_6_slope_fulldim_covers_4}{\href{\githubsearchurl?q=\%22def+kzh_6_slope_fulldim_covers_4(\%22}{\sage{kzh\underscore{}6\underscore{}slope\underscore{}fulldim\underscore{}covers\underscore{}4}}}
\pgfkeyssetvalue{/sagefunc/kzh_6_slope_fulldim_covers_5}{\href{\githubsearchurl?q=\%22def+kzh_6_slope_fulldim_covers_5(\%22}{\sage{kzh\underscore{}6\underscore{}slope\underscore{}fulldim\underscore{}covers\underscore{}5}}}
\pgfkeyssetvalue{/sagefunc/kzh_7_slope_1}{\href{\githubsearchurl?q=\%22def+kzh_7_slope_1(\%22}{\sage{kzh\underscore{}7\underscore{}slope\underscore{}1}}}
\pgfkeyssetvalue{/sagefunc/kzh_7_slope_2}{\href{\githubsearchurl?q=\%22def+kzh_7_slope_2(\%22}{\sage{kzh\underscore{}7\underscore{}slope\underscore{}2}}}
\pgfkeyssetvalue{/sagefunc/kzh_7_slope_3}{\href{\githubsearchurl?q=\%22def+kzh_7_slope_3(\%22}{\sage{kzh\underscore{}7\underscore{}slope\underscore{}3}}}
\pgfkeyssetvalue{/sagefunc/kzh_7_slope_4}{\href{\githubsearchurl?q=\%22def+kzh_7_slope_4(\%22}{\sage{kzh\underscore{}7\underscore{}slope\underscore{}4}}}
\pgfkeyssetvalue{/sagefunc/ll_strong_fractional}{\href{\githubsearchurl?q=\%22def+ll_strong_fractional(\%22}{\sage{ll\underscore{}strong\underscore{}fractional}}}
\pgfkeyssetvalue{/sagefunc/mlr_cpl3_a_2_slope}{\href{\githubsearchurl?q=\%22def+mlr_cpl3_a_2_slope(\%22}{\sage{mlr\underscore{}cpl3\underscore{}a\underscore{}2\underscore{}slope}}}
\pgfkeyssetvalue{/sagefunc/mlr_cpl3_b_3_slope}{\href{\githubsearchurl?q=\%22def+mlr_cpl3_b_3_slope(\%22}{\sage{mlr\underscore{}cpl3\underscore{}b\underscore{}3\underscore{}slope}}}
\pgfkeyssetvalue{/sagefunc/mlr_cpl3_c_3_slope}{\href{\githubsearchurl?q=\%22def+mlr_cpl3_c_3_slope(\%22}{\sage{mlr\underscore{}cpl3\underscore{}c\underscore{}3\underscore{}slope}}}
\pgfkeyssetvalue{/sagefunc/mlr_cpl3_d_3_slope}{\href{\githubsearchurl?q=\%22def+mlr_cpl3_d_3_slope(\%22}{\sage{mlr\underscore{}cpl3\underscore{}d\underscore{}3\underscore{}slope}}}
\pgfkeyssetvalue{/sagefunc/mlr_cpl3_f_2_or_3_slope}{\href{\githubsearchurl?q=\%22def+mlr_cpl3_f_2_or_3_slope(\%22}{\sage{mlr\underscore{}cpl3\underscore{}f\underscore{}2\underscore{}or\underscore{}3\underscore{}slope}}}
\pgfkeyssetvalue{/sagefunc/mlr_cpl3_g_3_slope}{\href{\githubsearchurl?q=\%22def+mlr_cpl3_g_3_slope(\%22}{\sage{mlr\underscore{}cpl3\underscore{}g\underscore{}3\underscore{}slope}}}
\pgfkeyssetvalue{/sagefunc/mlr_cpl3_h_2_slope}{\href{\githubsearchurl?q=\%22def+mlr_cpl3_h_2_slope(\%22}{\sage{mlr\underscore{}cpl3\underscore{}h\underscore{}2\underscore{}slope}}}
\pgfkeyssetvalue{/sagefunc/mlr_cpl3_k_2_slope}{\href{\githubsearchurl?q=\%22def+mlr_cpl3_k_2_slope(\%22}{\sage{mlr\underscore{}cpl3\underscore{}k\underscore{}2\underscore{}slope}}}
\pgfkeyssetvalue{/sagefunc/mlr_cpl3_l_2_slope}{\href{\githubsearchurl?q=\%22def+mlr_cpl3_l_2_slope(\%22}{\sage{mlr\underscore{}cpl3\underscore{}l\underscore{}2\underscore{}slope}}}
\pgfkeyssetvalue{/sagefunc/mlr_cpl3_n_3_slope}{\href{\githubsearchurl?q=\%22def+mlr_cpl3_n_3_slope(\%22}{\sage{mlr\underscore{}cpl3\underscore{}n\underscore{}3\underscore{}slope}}}
\pgfkeyssetvalue{/sagefunc/mlr_cpl3_o_2_slope}{\href{\githubsearchurl?q=\%22def+mlr_cpl3_o_2_slope(\%22}{\sage{mlr\underscore{}cpl3\underscore{}o\underscore{}2\underscore{}slope}}}
\pgfkeyssetvalue{/sagefunc/mlr_cpl3_p_2_slope}{\href{\githubsearchurl?q=\%22def+mlr_cpl3_p_2_slope(\%22}{\sage{mlr\underscore{}cpl3\underscore{}p\underscore{}2\underscore{}slope}}}
\pgfkeyssetvalue{/sagefunc/mlr_cpl3_q_2_slope}{\href{\githubsearchurl?q=\%22def+mlr_cpl3_q_2_slope(\%22}{\sage{mlr\underscore{}cpl3\underscore{}q\underscore{}2\underscore{}slope}}}
\pgfkeyssetvalue{/sagefunc/mlr_cpl3_r_2_slope}{\href{\githubsearchurl?q=\%22def+mlr_cpl3_r_2_slope(\%22}{\sage{mlr\underscore{}cpl3\underscore{}r\underscore{}2\underscore{}slope}}}
\pgfkeyssetvalue{/sagefunc/rlm_dpl1_extreme_3a}{\href{\githubsearchurl?q=\%22def+rlm_dpl1_extreme_3a(\%22}{\sage{rlm\underscore{}dpl1\underscore{}extreme\underscore{}3a}}}
\pgfkeyssetvalue{/sagefunc/automorphism}{\href{\githubsearchurl?q=\%22def+automorphism(\%22}{\sage{automorphism}}}
\pgfkeyssetvalue{/sagefunc/multiplicative_homomorphism}{\href{\githubsearchurl?q=\%22def+multiplicative_homomorphism(\%22}{\sage{multiplicative\underscore{}homomorphism}}}
\pgfkeyssetvalue{/sagefunc/projected_sequential_merge}{\href{\githubsearchurl?q=\%22def+projected_sequential_merge(\%22}{\sage{projected\underscore{}sequential\underscore{}merge}}}
\pgfkeyssetvalue{/sagefunc/restrict_to_finite_group}{\href{\githubsearchurl?q=\%22def+restrict_to_finite_group(\%22}{\sage{restrict\underscore{}to\underscore{}finite\underscore{}group}}}
\pgfkeyssetvalue{/sagefunc/interpolate_to_infinite_group}{\href{\githubsearchurl?q=\%22def+interpolate_to_infinite_group(\%22}{\sage{interpolate\underscore{}to\underscore{}infinite\underscore{}group}}}
\pgfkeyssetvalue{/sagefunc/two_slope_fill_in}{\href{\githubsearchurl?q=\%22def+two_slope_fill_in(\%22}{\sage{two\underscore{}slope\underscore{}fill\underscore{}in}}}
\pgfkeyssetvalue{/sagefunc/generate_example_e_for_psi_n}{\href{\githubsearchurl?q=\%22def+generate_example_e_for_psi_n(\%22}{\sage{generate\underscore{}example\underscore{}e\underscore{}for\underscore{}psi\underscore{}n}}}
\pgfkeyssetvalue{/sagefunc/chen_3_slope_not_extreme}{\href{\githubsearchurl?q=\%22def+chen_3_slope_not_extreme(\%22}{\sage{chen\underscore{}3\underscore{}slope\underscore{}not\underscore{}extreme}}}
\pgfkeyssetvalue{/sagefunc/psi_n_in_bccz_counterexample_construction}{\href{\githubsearchurl?q=\%22def+psi_n_in_bccz_counterexample_construction(\%22}{\sage{psi\underscore{}n\underscore{}in\underscore{}bccz\underscore{}counterexample\underscore{}construction}}}
\pgfkeyssetvalue{/sagefunc/gomory_fractional}{\href{\githubsearchurl?q=\%22def+gomory_fractional(\%22}{\sage{gomory\underscore{}fractional}}}
\pgfkeyssetvalue{/sagefunc/not_minimal_2}{\href{\githubsearchurl?q=\%22def+not_minimal_2(\%22}{\sage{not\underscore{}minimal\underscore{}2}}}
\pgfkeyssetvalue{/sagefunc/not_extreme_1}{\href{\githubsearchurl?q=\%22def+not_extreme_1(\%22}{\sage{not\underscore{}extreme\underscore{}1}}}
\pgfkeyssetvalue{/sagefunc/kzh_2q_example_1}{\href{\githubsearchurl?q=\%22def+kzh_2q_example_1(\%22}{\sage{kzh\underscore{}2q\underscore{}example\underscore{}1}}}
\pgfkeyssetvalue{/sagefunc/extremality_test}{\href{\githubsearchurl?q=\%22def+extremality_test(\%22}{\sage{extremality\underscore{}test}}}
\pgfkeyssetvalue{/sagefunc/plot_2d_diagram}{\href{\githubsearchurl?q=\%22def+plot_2d_diagram(\%22}{\sage{plot\underscore{}2d\underscore{}diagram}}}
\pgfkeyssetvalue{/sagefunc/nice_field_values}{\href{\githubsearchurl?q=\%22def+nice_field_values(\%22}{\sage{nice\underscore{}field\underscore{}values}}}
\pgfkeyssetvalue{/sagefunc/ParametricRealFieldElement}{\href{\githubsearchurl?q=\%22def+ParametricRealFieldElement(\%22}{\sage{ParametricRealFieldElement}}}
\pgfkeyssetvalue{/sagefunc/ParametricRealField}{\href{\githubsearchurl?q=\%22def+ParametricRealField(\%22}{\sage{ParametricRealField}}}

%% file: cpl_extreme_cases_diagrams.tex
\begin{figure}[h]
\centering
\includegraphics[width=0.5\textwidth]{\groupcodedir/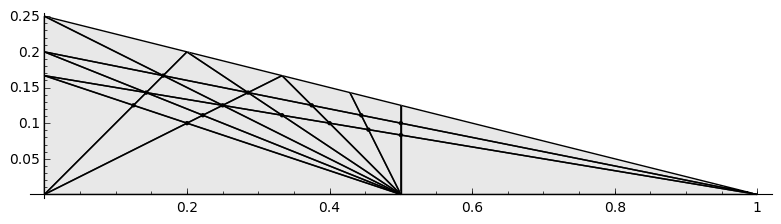}
\caption{The plane of the parameters $(r_0, z_1)$, which determine the
  location of breakpoints, has 30 two-dimensional cells derived from the
  arrangement of breakpoints of the CPL function.  Together with the
  lower-dimensional cells that arise as intersections, there are a total of 87
  cells.} 
\label{fig:region_87}
\end{figure}

\begin{figure}[h]
\begin{minipage}[t]{.49\textwidth}
\includegraphics[width=\linewidth]{\groupcodedir/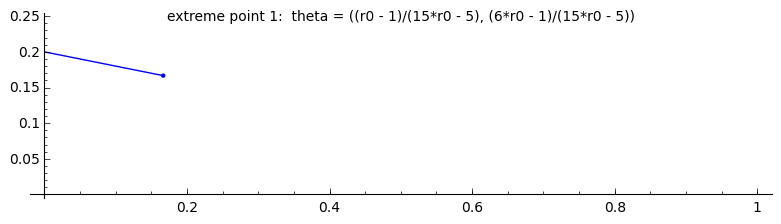}
\center{$i=1$, \cite[Table 5]{Miller-Li-Richard2008} ext pt `f'}
\includegraphics[width=\linewidth]{\groupcodedir/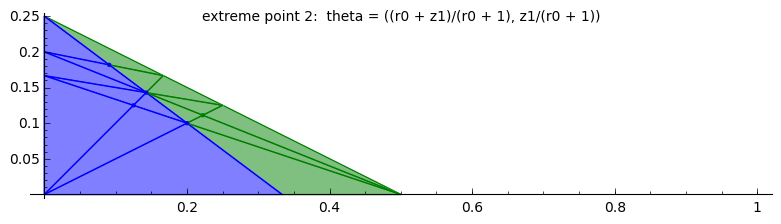}
\center{$i=2$, \cite[Table 5]{Miller-Li-Richard2008} ext pt `c' }
\includegraphics[width=\linewidth]{\groupcodedir/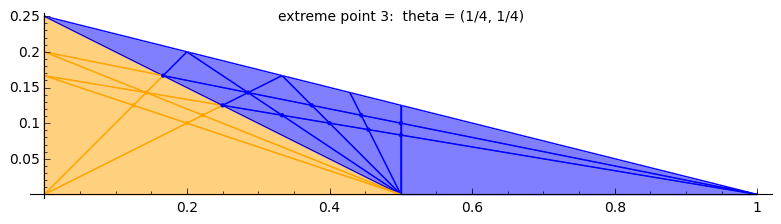}
\center{$i=3$, \cite[Table 5]{Miller-Li-Richard2008} ext pt `h' }
\includegraphics[width=\linewidth]{\groupcodedir/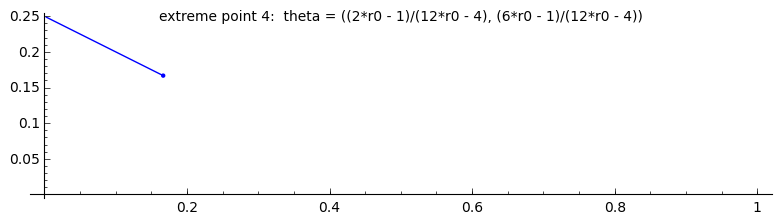}
\center{$i=4$, \cite[Table 5]{Miller-Li-Richard2008} ext pt `g' }
\includegraphics[width=\linewidth]{\groupcodedir/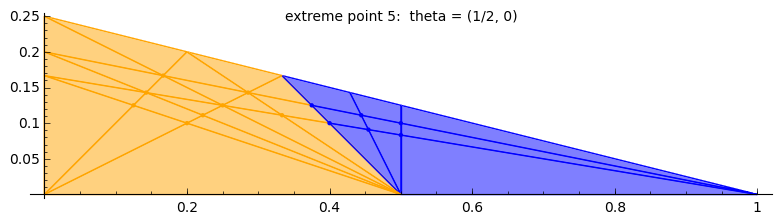}
\center{$i=5$, \cite[Table 5]{Miller-Li-Richard2008} ext pt `r, o' }
\end{minipage}
\begin{minipage}[t]{.49\textwidth}
\includegraphics[width=\linewidth]{\groupcodedir/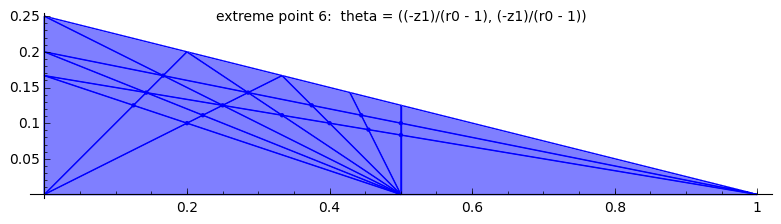}
\center{$i=6$, \cite[Table 5]{Miller-Li-Richard2008} ext pt `a' }
\includegraphics[width=\linewidth]{\groupcodedir/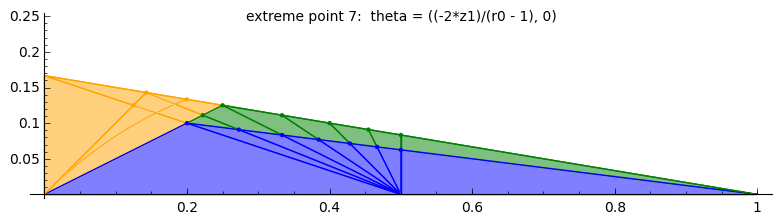}
\center{$i=7$, \cite[Table 5]{Miller-Li-Richard2008} ext pt `n, d' }
\includegraphics[width=\linewidth]{\groupcodedir/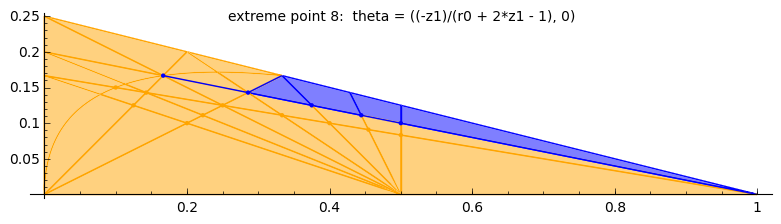}
\center{$i=8$, \cite[Table 5]{Miller-Li-Richard2008} ext pt `q, p, l, k'}
\includegraphics[width=\linewidth]{\groupcodedir/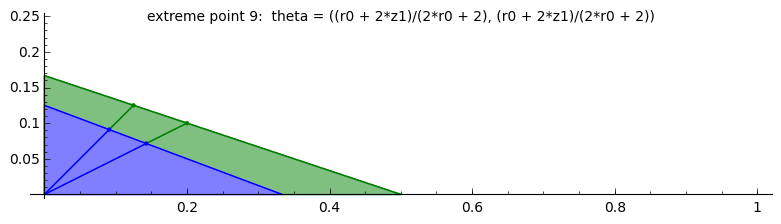}
\center{$i=9$, \cite[Table 5]{Miller-Li-Richard2008} ext pt `b' }
\end{minipage}
\caption{The plane of the parameters $(r_0, z_1)$, which determine the
  location of breakpoints, with cells corresponding to the automatic proof. 
  Cells are color-coded: `not constructible' (\emph{white}), `constructible,
  but not minimal' (\emph{yellow}), `minimal, but not
  extreme' (\emph{green}) or `extreme' (\emph{blue})}
\label{fig:cpl_theta_i}
\end{figure}

%% file: param_3s_ext_1_f_diagrams.tex
\begin{figure}[h]
\centering
\begin{minipage}[b]{.19\textwidth}
\includegraphics[width=\linewidth]{\groupcodedir/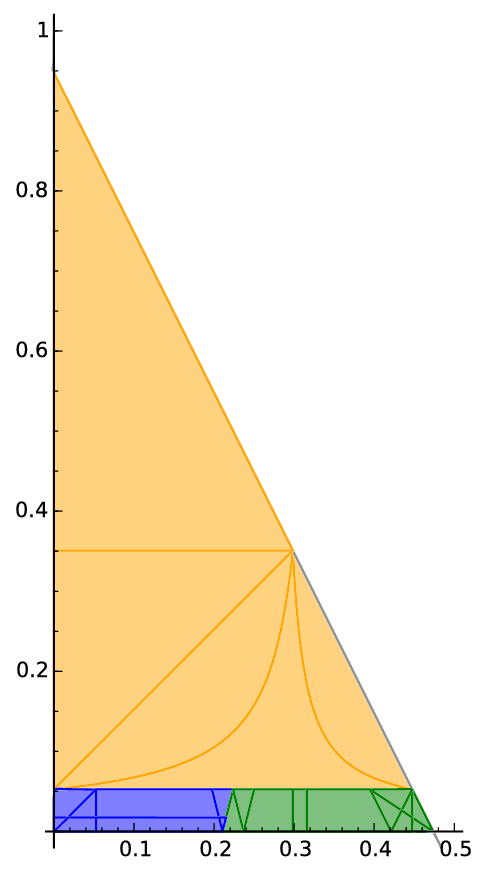}
\center{$i=1$}
\includegraphics[width=\linewidth]{\groupcodedir/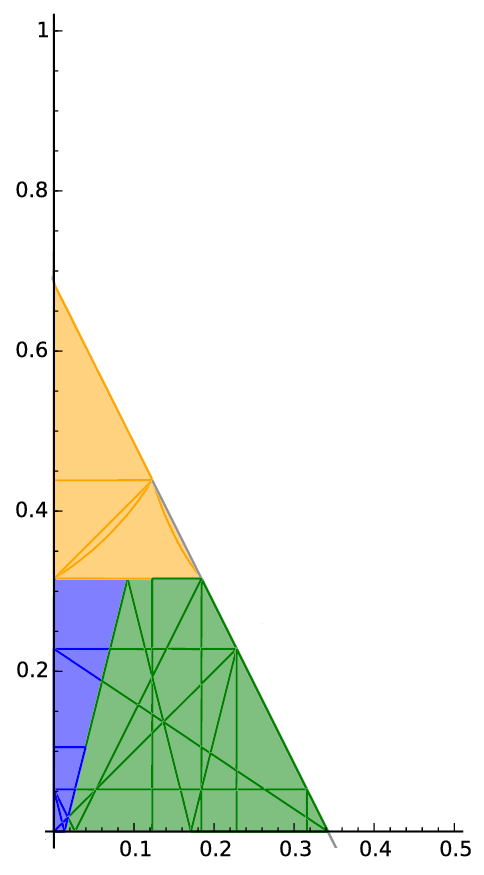}
\center{$i=6$}
\end{minipage}
\begin{minipage}[b]{.19\textwidth}
\includegraphics[width=\linewidth]{\groupcodedir/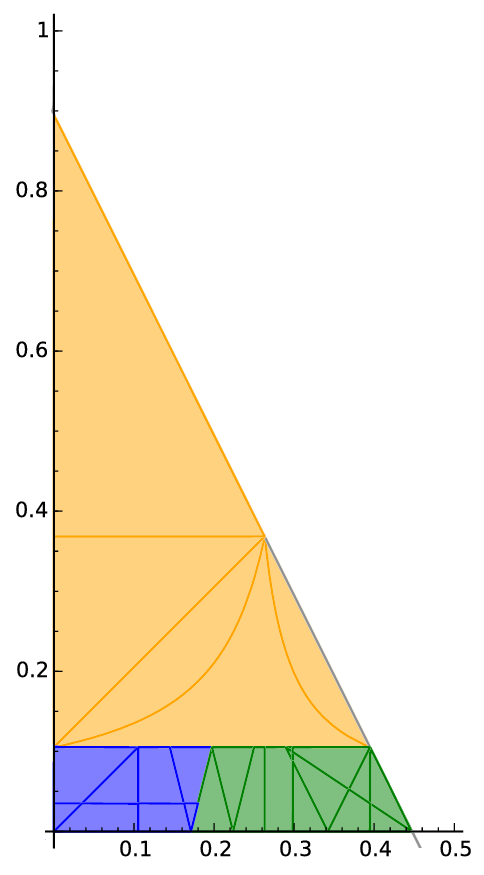}
\center{$i=2$}
\includegraphics[width=\linewidth]{\groupcodedir/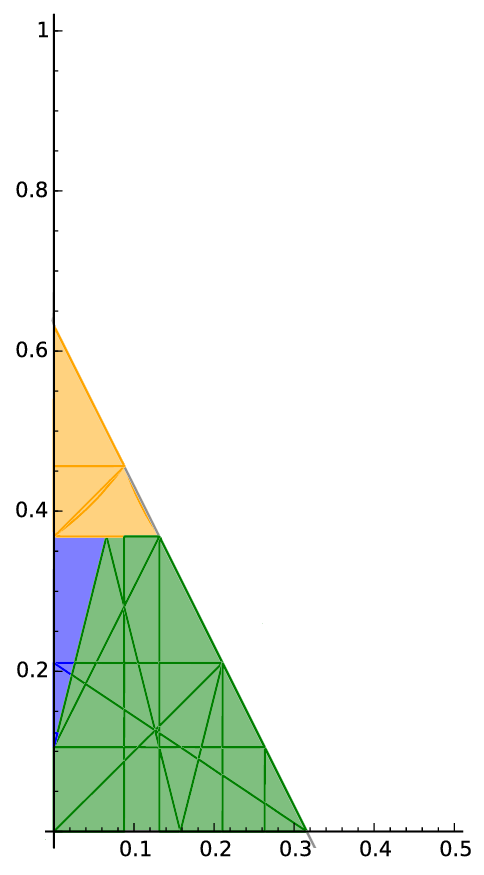}
\center{$i=7$}
\end{minipage}
\begin{minipage}[b]{.19\textwidth}
\includegraphics[width=\linewidth]{\groupcodedir/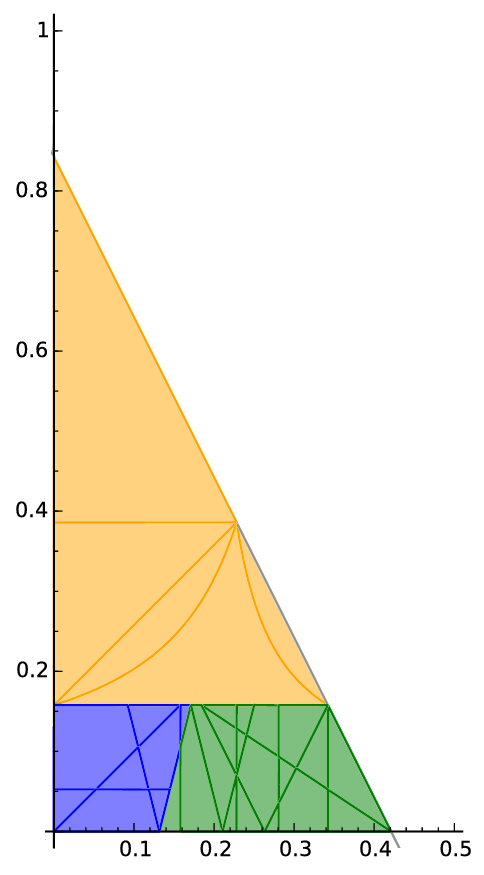}
\center{$i=3$}
\includegraphics[width=\linewidth]{\groupcodedir/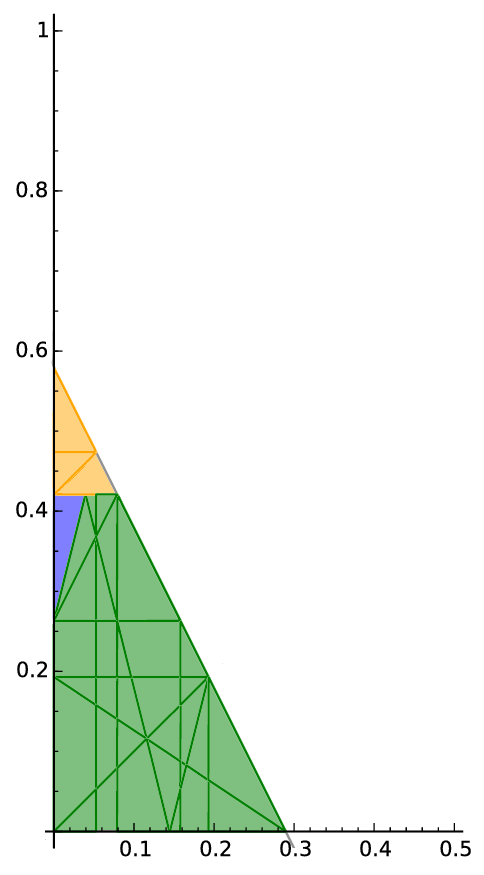}
\center{$i=8$}
\end{minipage}
\begin{minipage}[b]{.19\textwidth}
\includegraphics[width=\linewidth]{\groupcodedir/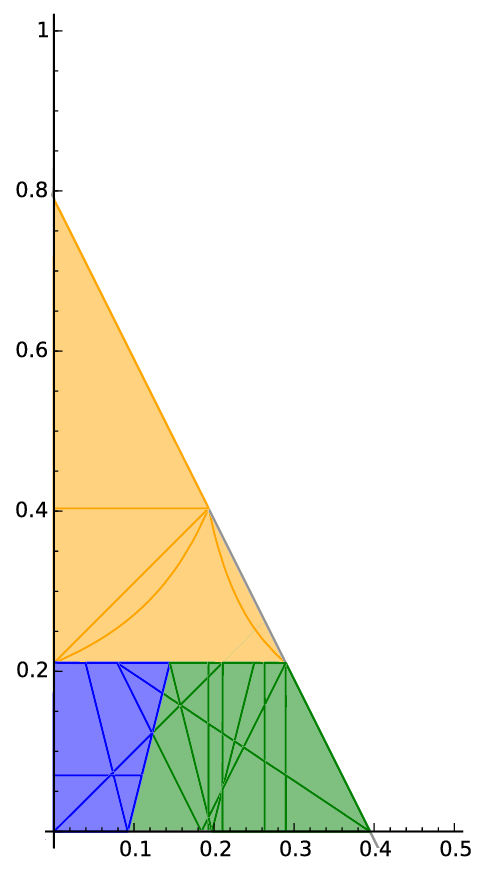}
\center{$i=4$}
\includegraphics[width=\linewidth]{\groupcodedir/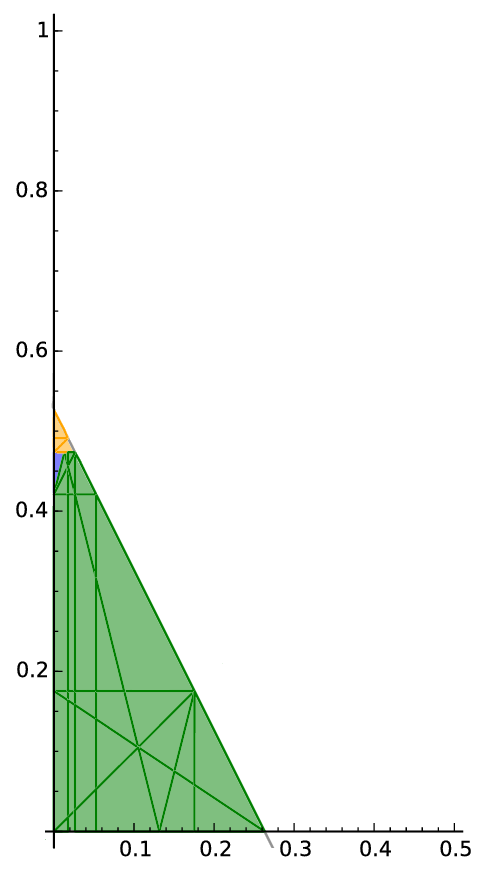}
\center{$i=9$}
\end{minipage}
\begin{minipage}[b]{.19\textwidth}
\includegraphics[width=\linewidth]{\groupcodedir/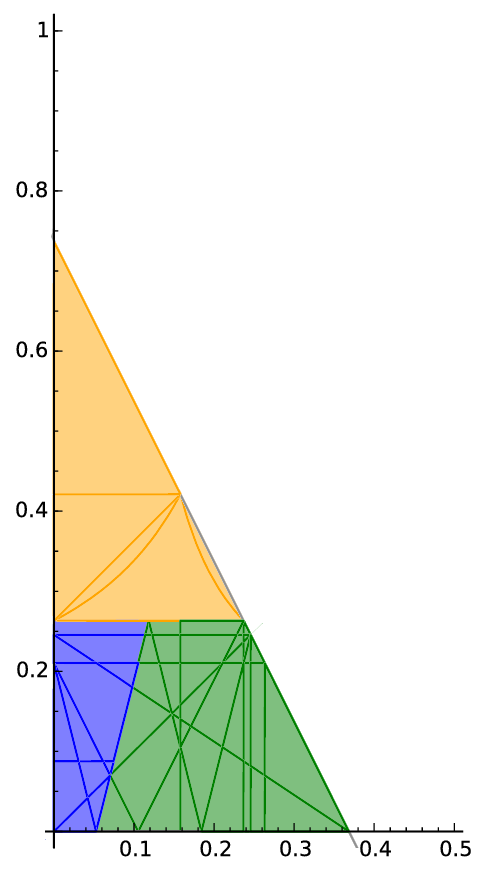}
\center{$i=5$}
\includegraphics[width=\linewidth]{\groupcodedir/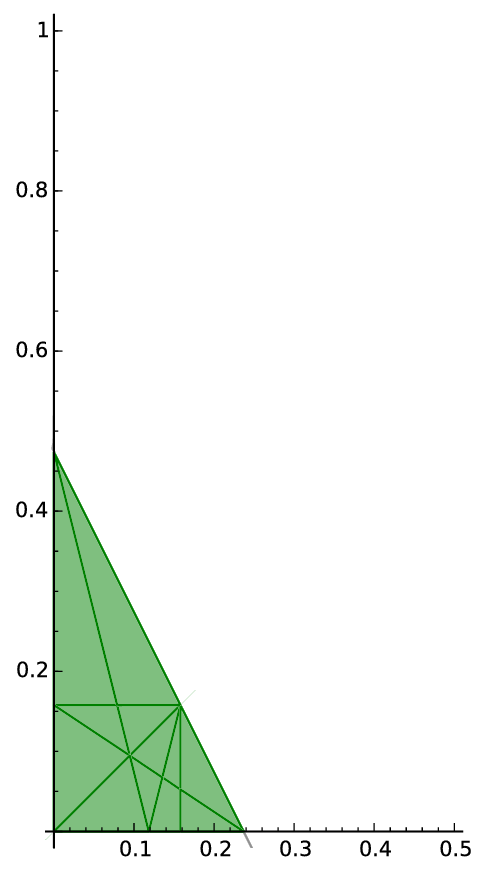}
\center{$i=10$}
\end{minipage}
\caption{Slices of the cell complex for the family \sage{kzh\_3\_slope\_param\_extreme\_1} for $f = i/19$}
\label{fig:kzh_3_slope_param_extreme_1_slices}
\end{figure}


%% file: param_3s_ext_2_f_diagrams.tex
\begin{figure}[h]
\centering
\begin{minipage}[b]{.32\textwidth}
\includegraphics[width=\linewidth]{\groupcodedir/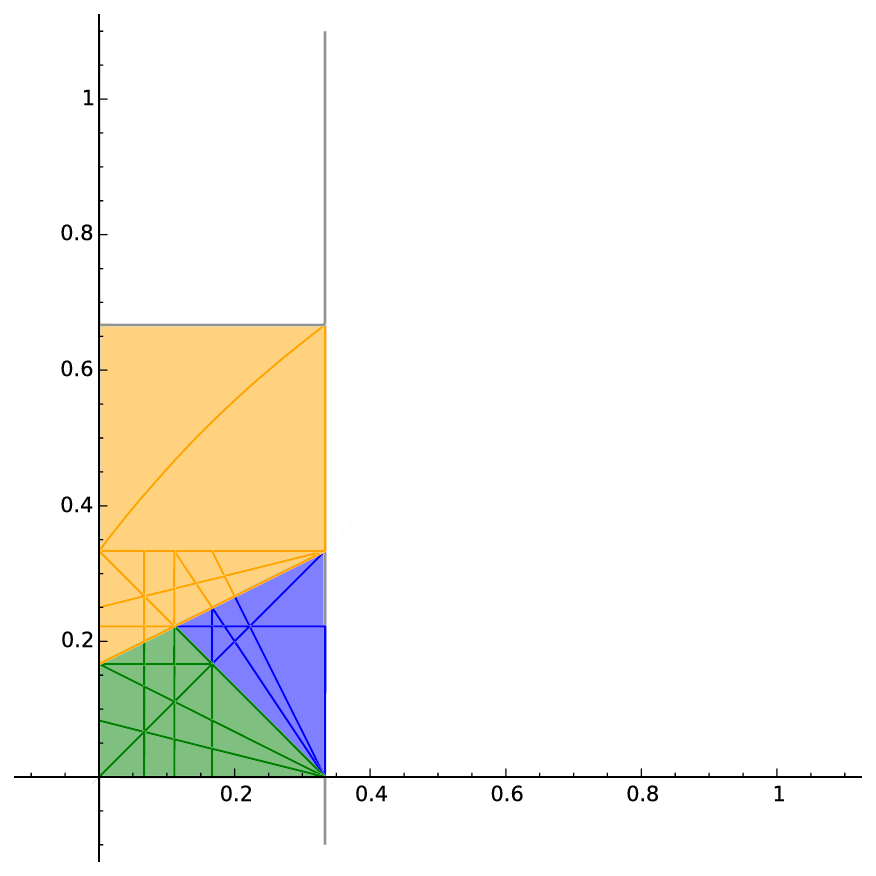}
\center{$i=3$}
\includegraphics[width=\linewidth]{\groupcodedir/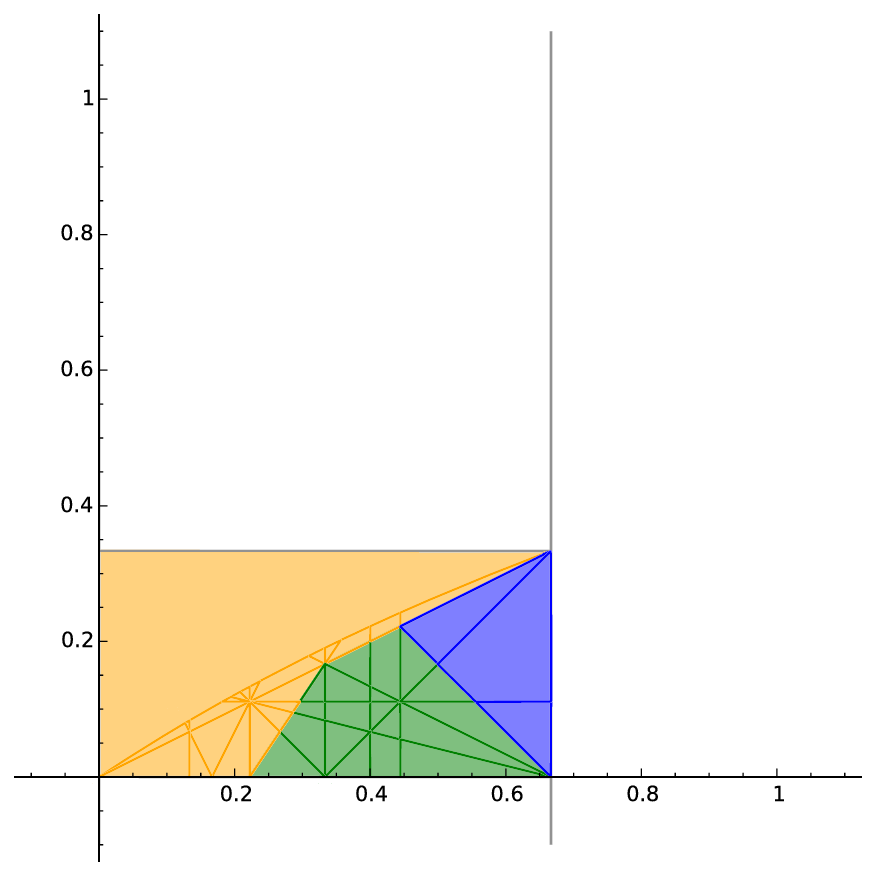}
\center{$i=6$}
\end{minipage}
\begin{minipage}[b]{.32\textwidth}
\includegraphics[width=\linewidth]{\groupcodedir/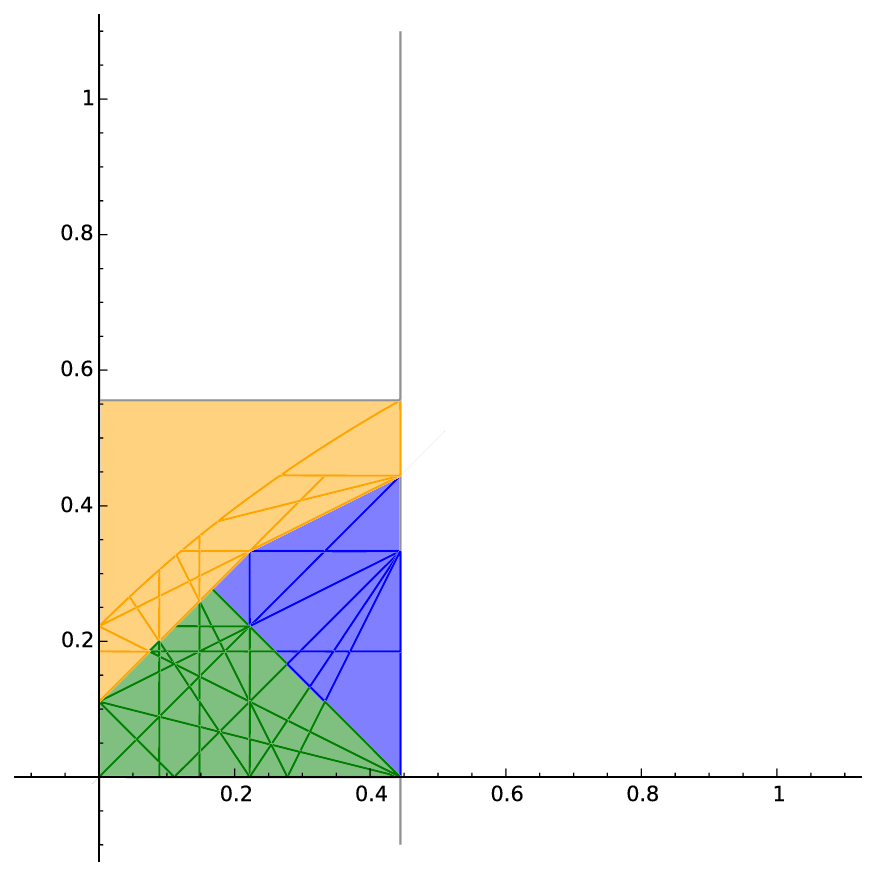}
\center{$i=4$}
\includegraphics[width=\linewidth]{\groupcodedir/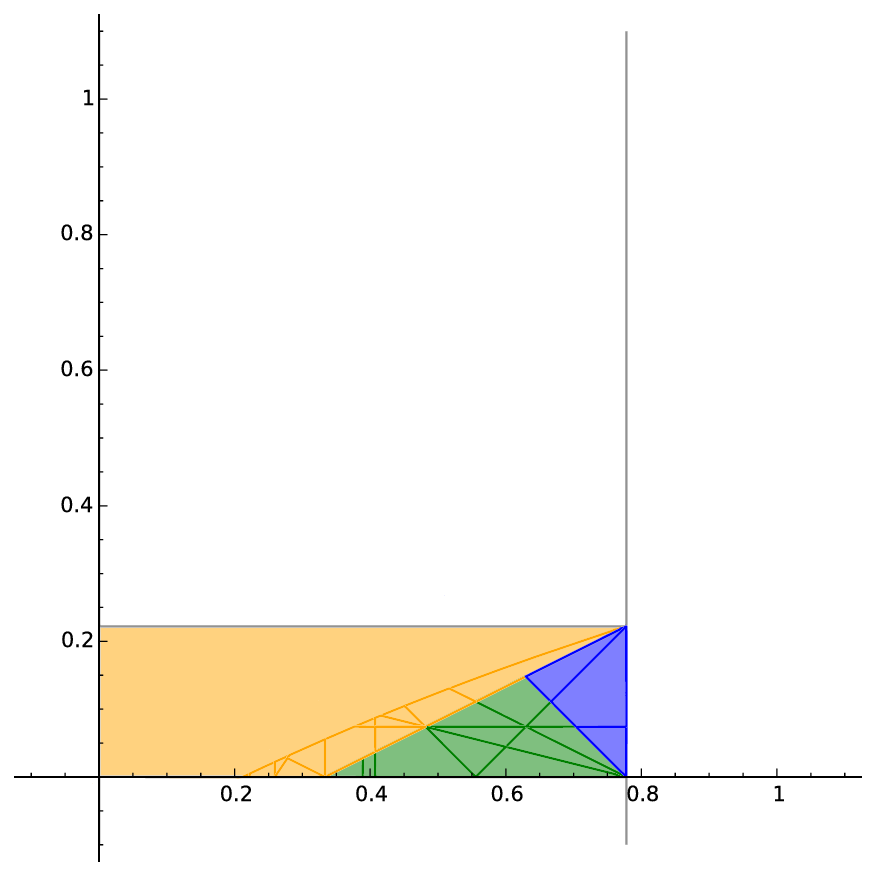}
\center{$i=7$}
\end{minipage}
\begin{minipage}[b]{.32\textwidth}
\includegraphics[width=\linewidth]{\groupcodedir/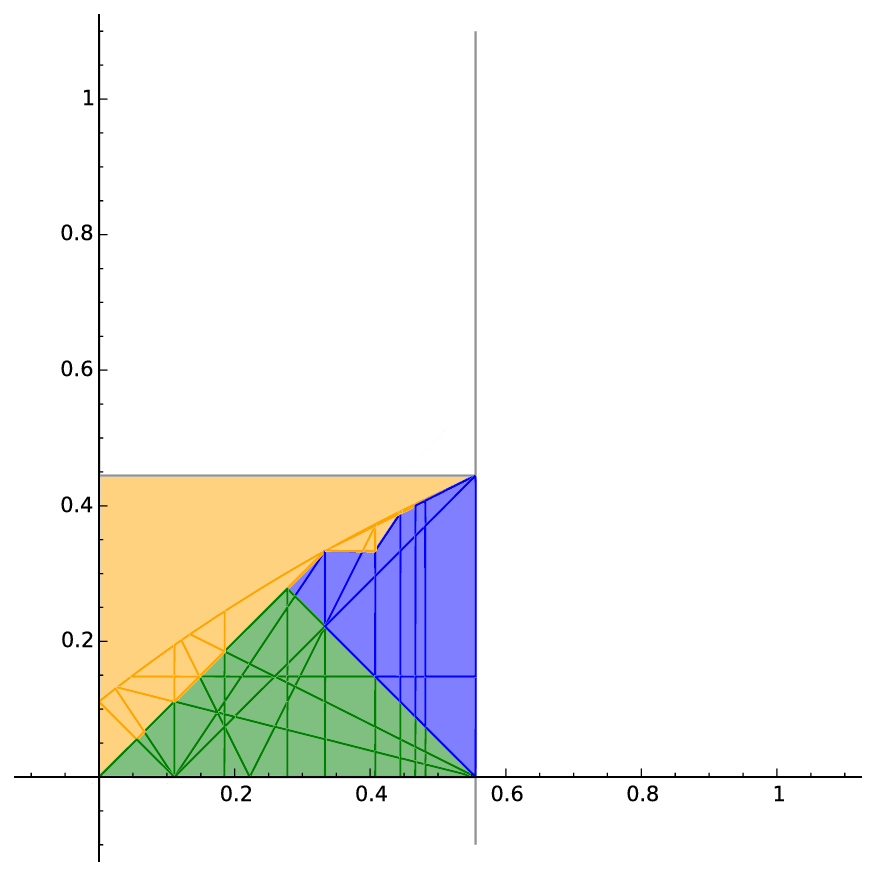}
\center{$i=5$}
\includegraphics[width=\linewidth]{\groupcodedir/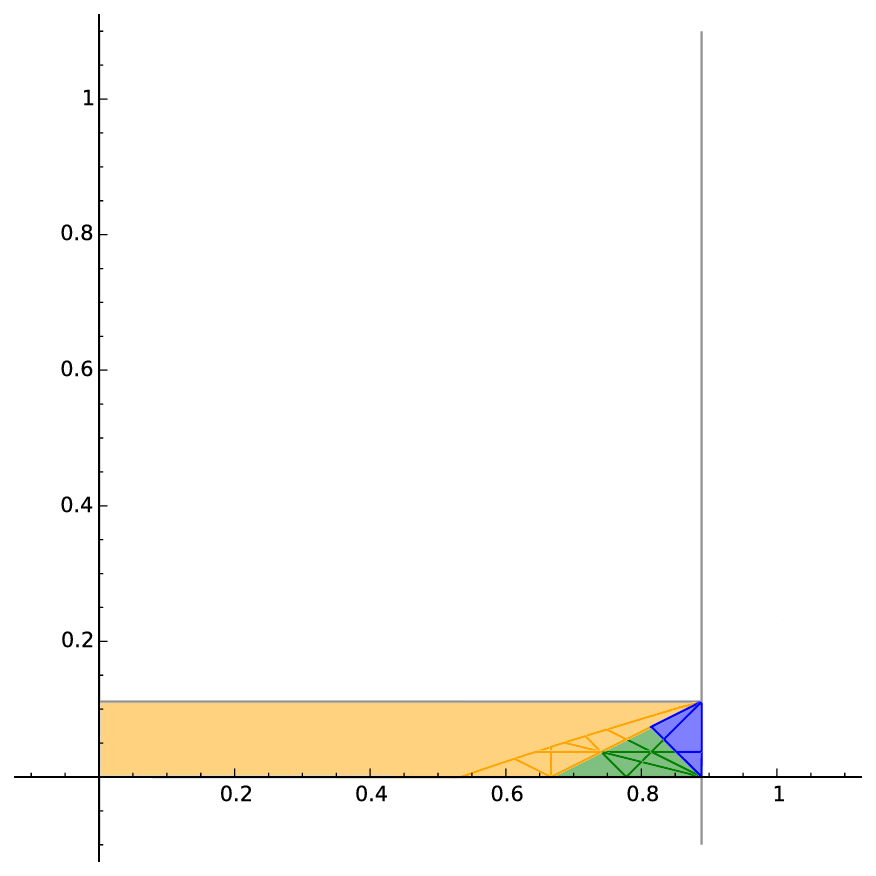}
\center{$i=8$}
\end{minipage}
\caption{Slices of the cell complex for the family \sage{kzh\_3\_slope\_param\_extreme\_2} for $f = i/9$}
\label{fig:kzh_3_slope_param_extreme_2_slices}
\end{figure}